\documentclass[12pt]{amsart}
\usepackage{amsfonts}
\usepackage{amssymb}
\newtheorem{thm}{Theorem}[section]
\newtheorem{cor}[thm]{Corollary}
\newtheorem{lem}[thm]{Lemma}
\newtheorem{prop}[thm]{Proposition}
\theoremstyle{definition}
\newtheorem{defn}[thm]{Definition}
\theoremstyle{remark}
\newtheorem{rem}[thm]{Remark}
\newcommand{\R}{\mathbb R}
\newcommand{\C}{\mathbb C}

\newcommand{\To}{\Longrightarrow}
\newcommand{\X}{\mathfrak X}

\begin{document}

\title[K\"AHLER MANIFOLDS OF QUASI-CONSTANT HOLOMORPHIC CURVATURES]
{K\"AHLER MANIFOLDS OF QUASI-CONSTANT HOLOMORPHIC SECTIONAL CURVATURES}%
\author{G. Ganchev and V. Mihova}%
\address{Bulgarian Academy of Sciences, Institute of Mathematics and Informatics,
Acad. G. Bonchev Str. bl. 8, 1113 Sofia, Bulgaria}%
\email{ganchev@math.bas.bg}%
\address{Faculty of Mathematics and Informatics, University of Sofia,
J. Bouchier Str. 5, (1164) Sofia, Bulgaria}
\email{mihova@fmi.uni-sofia.bg}

\begin{abstract}
The K\"ahler manifolds of quasi-constant holomorphic sectional
curvatures are introduced as K\"ahler manifolds with complex
distribution of codimension two, whose holomorphic sectional
curvature only depends on the corresponding point and the geometric
angle, associated with the section. A curvature identity
characterizing such manifolds is found. The biconformal group of
transformations whose elements transform K\"ahler metrics into
K\"ahler ones is introduced and biconformal tensor invariants are
obtained. This makes it possible to classify the manifolds under
consideration locally. The class of locally biconformal flat K\"ahler
metrics is shown to be exactly the class of K\"ahler metrics whose
potential function is only a function of the distance from the origin
in ${\C}^n$. Finally we show that any rotational even dimensional
hypersurface carries locally a natural K\"ahler structure which is
of quasi-constant holomorphic sectional curvatures.
\end{abstract}
\maketitle
\section{Introduction}

In the present paper we study K\"ahler manifolds $(M,g,J,D)\;(\dim M = 2n)$
with $J$-invariant distribution $D$ of $\dim D = 2(n-1)$. The orthogonal
distribution $D^{\perp}$ of the given distribution $D$ is also $J$-invariant
and $\dim D^{\perp}=2$. The structural group of these manifolds is
$U(n-1)\times U(1)$.

For example any real function $u \in C^{\infty}$ on $M$ with $du
\neq 0$ generates a $J$-invariant distribution $D$ with
$D^{\perp}=span\{grad \, u, Jgrad \, u\}$. Thus any K\"ahler
manifold with scalar curvature $\tau$ in general carries such a
structure $(D,D^{\perp})$ generated by $grad \, \tau$ and $Jgrad
\, \tau$. In this case the structure $(D,D^{\perp})$ is determined
by the K\"ahler structure $(g,J)$.

In Section 2 we introduce the basic notion of a K\"ahler manifold
$(M,g,J,D)$ $(n\geq2)$ with quasi-constant holomorphic sectional
curvatures (K\"ahler $QCH$-manifolds).

Any holomorphic tangent section at a point $p \in M$ and
$D^{\perp}(p)$ determines a geometric angle $\varphi$. Then the
notion of a K\"ahler manifold of quasi-constant holomorphic
sectional curvatures is introduced by the natural condition:
\vskip 1mm
{\it The holomorphic sectional curvatures of the manifold $(M,g,J,D)$ depend only
on the angle $\varphi$ and the point $p \in M$.}
\vskip 1mm
The notion of a K\"ahler $QCH$-manifold is the K\"ahler analogue of the notion
of a Riemannian manifold of quasi-constant sectional curvatures \cite {BP, GM}.

We construct three K\"ahler tensors $\pi, \Phi$ and $\Psi$ invariant under the
action of the structural group and prove Proposition \ref {P:2.2}, which gives
the main tool to investigate the K\"ahler $QCH$-manifolds:
\vskip 1mm
{\it A K\"ahler manifold $(M,g,J,D)$ is of quasi-constant holomorphic sectional
curvatures if and only if its curvature tensor $R$ satisfies the identity
$$R=a\pi+b\Phi+c\Psi,\leqno{(1.1)}$$
where $a, b$ and $c$ are functions on $M$.}

In Section 3 we study the integrability conditions for the
curvature identity (1.1). A complete system of integrability
conditions for (1.1) is obtained in Theorem \ref{T:3.4}. Applying
this theorem we obtain:
\vskip 1mm
{\it If the distribution $D$ of a K\"ahler $QCH$-manifold is not involutive,
then the integrability conditions of $(1.1)$ reduce to the conditions
$(3.21)$.}
\vskip 1mm
In Section 4 we discuss transformations of K\"ahler metrics into
K\"ahler ones. By using the induced metric $\eta \otimes \eta +
\tilde\eta\otimes\tilde\eta$ on $D^{\perp}$, we introduce biconformal
transformations of the structure $(g, \eta)$ by the formulas
$$g'=e^{2u} \{g+(e^{2v}-1)(\eta \otimes \eta + \tilde\eta\otimes\tilde\eta)\},
\quad \eta'=e^{u+v}\eta,\leqno{(1.2)}$$
where $u$ and $v$ are functions on $M$ satisfying certain conditions.

We introduce the notion of a $B$-distribution with the conditions (4.15) and
prove that:
\vskip 1mm
{\it All metrics $g'$ given by $(1.2)$ are K\"ahlerian if and only if the
distribution $D$ is a $B$-distribution.}
\vskip 1mm

Further we define $B_0$-distributions as a special case of $B$-distributions with
the conditions (4.16).
In Proposition \ref{P:4.3} we show that the transformations (1.2) form a group -
the group of biconformal transformations.

Section 5 is devoted to the tensor invariants of the biconformal
group in the case of a $B_0$-distribution $D$. Theorem \ref{T:5.5}
states as follows:
\vskip 1mm
{\it The tensor $R-a\pi-b\Phi-c\Psi$ of type $(1,3)$ is a biconformal invariant.}
\vskip 1mm
This theorem makes possible a local classification of the K\"ahler $QCH$-manifolds
satisfying a certain inequality (Theorem \ref{T:5.8}):
\vskip 1mm
{\it A K\"ahler manifold $(M,g,J,D)$ with $B_0$-distribution $D$ is biconformally
flat if and only if}
$$R=a\pi +b\Phi +c\Psi, \quad a+k^2>0.$$
\vskip 1mm
The main property, which connects the K\"ahler $QCH$-manifolds with the set of
all K\"ahler metrics, whose potential function is of the type
$f(r^2)$\;($r$ - the distance from the origin in ${\C}^n$) is established by
Theorem \ref{T:5.10}:
\vskip 1mm
{\it Any K\"ahler metric $g=\partial\bar\partial f(r^2)$ is biconformally
flat and vice versa.}
\vskip 1mm
In Section 6 we show that any rotational hypersurface $(M^{2n},\bar g)$ in
${\C}^n\times {\R}$ with axis of revolution $l={\R}$, which has no common
points with $l$ carries a geometrically determined complex structure $J$.
Thus $(M^{2n},\bar g, J)$ can be considered as a locally conformal K\"ahler
manifold in all dimensions $2n\geq 4$.

Further we introduce a natural K\"ahler metric on $(M^{2n},\bar g, J)$ by
the formula (6.16). This formula makes the class of rotational hypersurfaces
an important source of K\"ahler metrics because of Theorem \ref {T:6.2}:
\vskip 1mm
{\it Let $(M^{2n},\bar g,J,\bar\xi) \, (2n\geq 4)$ be a rotational hypersurface
satisfying the conditions $(6.15)$. Then the complex dilatational K\"ahler metric
$g$ given by $(6.16)$ is of quasi-constant holomorphic sectional curvatures.}
\vskip 1mm
Finally we find the rotational hypersurfaces $M^{2n}$ whose complex dilatational
K\"ahler metric is of constant holomorphic sectional curvatures:
\vskip 1mm
{\it Any rotational hypersurface $M^{2n}$ which carries a complex dilatational
K\"ahler metric of constant holomorphic sectional curvature $a=const>0$ is
generated by a meridian of the type
$$\gamma : y=\pm\frac{1}{\sqrt a}\left(\sqrt{8-ax^2}+\ln{\frac{\sqrt{8-ax^2}-
2}{\sqrt{8-ax^2}+2}}\right)+y_0, \quad 0<x<\frac{2}{\sqrt a},$$
where the meridian $\gamma$ is considered with respect to the usual coordinate
system $Oxy$ with axis of revolution $l=Oy$.}

\vskip 2mm
\section{A tensor characterization of the K\"ahler manifolds of quasi-constant
         holomorphic sectional curvatures}
\vskip 2mm

Let $(M,g,J,D)$ be a 2n-dimensional K\"ahler manifold with metric $g$,
complex structure $J$ and $J$-invariant distribution $D$ of codimension 2.
The Lie algebra of all $C^{\infty}$ vector fields on $M$ will be
denoted by ${\X}M$ and $T_pM $ will stand for the tangent space to
$M$ at any point $p \in M$. In the presence of the distribution
$D$ the structure of any tangent space is $T_pM = D(p) \oplus
D^{\perp}(p)$, where $D^{\perp}(p)$ is the 2-dimensional
$J$-invariant orthogonal complement to the space $D(p)$. This
means that the structural group of the manifolds under
consideration is the subgroup $U(n-1)\times U(1)$ of $U(n)$.

As our considerations are local, we can assume the existence of a
unit vector field $\xi$ on $M$ such that $D^{\perp}(p) = span \{\xi, J\xi\}$
at any point $p \in M $. We denote by $\eta$ and
$\tilde \eta$ the unit 1-forms corresponding to $\xi$ and $J\xi$,
respectively, i.e.
$$\eta(X) = g(\xi, X), \quad \tilde \eta (X) = g(J\xi, X) = -\eta(JX);
\quad X \in {\X}M.$$
Then the distribution $D$ is determined by the conditions
$$D(p) = \{X \in T_pM \, \vert \, \eta(X) = \tilde \eta(X) = 0 \}, \quad
p \in M.$$

As a rule, we use the following denotations for vector fields (vectors):
$$X,Y,Z \in {\X}M \; (T_pM); \quad x_0,y_0,z_0 \in {\X}D \; (D(p)).$$

The K\"ahler form $\Omega$ of the structure $(g,J)$ is given by
$\Omega(X,Y) = g(JX,Y)$, $X,Y \in {\X}M$.

Let $\nabla$ be the Levi-Civita connection of the metric $g$. The
Riemannian curvature tensor $R$, the Ricci tensor $\rho$ and the
scalar curvature $\tau$ of $\nabla$ are given by the equalities
$$R(X,Y)Z = \nabla_X\nabla_YZ - \nabla_Y\nabla_XZ - \nabla_{[X,Y]}Z,$$
$$R(X,Y,Z,U) = g(R(X,Y)Z,U); \quad X,Y,Z,U \in {\X}M,$$
$$\rho (Y,Z) = \displaystyle{\sum_{i=1}^{2n}} R(e_i, Y, Z, e_i);
\quad Y,  Z \in T_pM,$$
$$\tau = \displaystyle{\sum_{i = 1}^{2n}\rho(e_i,e_i)},$$
where $\{e_i\}, \, i = 1,...,2n$ is an orthonormal basis for $T_pM, \, p \in M.$

We recall that the curvature tensor $R$ of any K\"ahler manifold
satisfies the identities
$$\begin{array}{l}
R(X,Y,Z,U) = - R(Y,X,Z,U);\\[2mm]
R(X,Y)Z + R(Y,Z)X + R(Z,X)Y = 0;\\[2mm]
R(X,Y,Z,U) = - R(X,Y,U,Z);\\[2mm]
R(X,Y)JZ = JR(X,Y)Z.
\end{array} \leqno{(2.1)}$$

Now we shall introduce the geometric functions and tensors
associated with the structures $(g,J,D)$.

All directions in $span\{\xi, J\xi\}$ have one and the same Ricci
curvature which is denoted by $\sigma$, i. e.
$$\sigma = \rho (\xi, \xi) = \rho (J\xi, J\xi).\leqno{(2.2)}$$

The Riemannian sectional curvature of the structural distribution
$D^{\perp}$ is denoted by $\varkappa$, i. e.
$$\varkappa = R(\xi, J\xi, J\xi, \xi).\leqno{(2.3)}$$

Thus the structures $(g,J,D)$ give rise to the functions $\varkappa,
\sigma$ and $\tau$.

Further we note that the tensor $\eta(X)\eta(Y) +
\tilde\eta(X)\tilde\eta(Y)$ does not depend on the basis
$\{\xi,J\xi\}$. Then the fundamental symmetric tensors of type
(0,2) are
$$g(X,Y), \quad \eta(X)\eta(Y) + \tilde \eta(X)
\tilde \eta(Y); \quad X,Y \in {\X}M.$$

We also mention the corresponding fundamental skew symmetric tensors
$$\Omega (X,Y), \quad \eta(X)\tilde\eta(Y)
- \eta(Y)\tilde\eta(X); \quad X,Y \in {\X}M.$$

Any tensor over $T_pM, \, p \in M$ of type (0,4) having the
symmetries (2.1) is called a K\"ahler tensor. We need the
following invariant K\"ahler tensors:
$$\begin{array}{ll}
4\pi (X,Y,Z,U) = &g(Y,Z)g(X,U) - g(X,Z)g(Y,U)\\[2mm]
&+ g(JY,Z)g(JX,U) - g(JX,Z)g(JY,U)\\[2mm]
&- 2g(JX,Y)g(JZ,U);
\end{array}\leqno{(2.4)}$$
$$\begin{array}{ll}
8\Phi (X,Y,Z,U) =
&g(Y,Z)\{\eta(X)\eta(U) + \tilde\eta(X)\tilde\eta(U)\}\\
[2mm]
&- g(X,Z)\{\eta(Y)\eta(U)+ \tilde\eta(Y)\tilde\eta(U)\}\\
[2mm]
&+ g(X,U)\{\eta(Y)\eta(Z) + \tilde\eta(Y)\tilde\eta(Z)\}\\
[2mm]
& - g(Y,U)\{\eta(X)\eta(Z) + \tilde\eta(X)\tilde\eta(Z)\}\\
[2mm]
& + g(JY,Z)\{\eta(X)\tilde\eta(U) - \eta(U)\tilde\eta(X)\}\\
[2mm]
&- g(JX,Z)\{\eta(Y)\tilde\eta(U) - \eta(U)\tilde\eta(Y)\}\\
[2mm]
& + g(JX,U)\{\eta(Y)\tilde\eta(Z) - \eta(Z)\tilde\eta(Y)\}\\
[2mm]
&- g(JY,U)\{\eta(X)\tilde\eta(Z) - \eta(Z)\tilde\eta(X)\}\\
[2mm]
&- 2g(JX,Y)\{\eta(Z)\tilde\eta(U) - \eta(U)\tilde\eta(Z)\}\\
[2mm] &- 2g(JZ,U)\{\eta(X)\tilde\eta(Y) -
\eta(Y)\tilde\eta(X)\};$$
\end{array}\leqno{(2.5)}$$
$$\begin{array}{lll}
\Psi (X,Y,Z,U)& = &\eta(Y)\eta(Z)\tilde\eta(X)\tilde\eta(U)
- \eta(X)\eta(Z)\tilde\eta(Y)\tilde\eta(U)\\
[2mm] & &+ \eta(X)\eta(U)\tilde\eta(Y)\tilde\eta(Z) -
\eta(Y)\eta(U)\tilde\eta(X)\tilde\eta(Z)\\
[2mm] &=&\{(\eta \wedge
\tilde\eta)\otimes(\tilde\eta\wedge\eta)\}(X,Y,Z,U),
\end{array}\leqno{(2.6)}$$
$X,Y,Z,U \in {\X}M.$

These tensors are invariant under the action of the structural
group $U(n-1)\times U(1)$ on the tensors over $T_pM, \, p \in M$
in the standard sense (e.g. \cite{TV}).

Let $(M,g,J,D) \,(\dim \, M = 2n)$ be a K\"ahler manifold with J-invariant
distributions $D \,(\dim \, D = 2(n-1))$ and $D^{\perp}$. The structures
$(g,J,D)$ give rise to the following {\it geometric} types of sectional
curvatures with respect to any K\"ahler curvature tensor (e. g. $R$):

{\it horizontal} sectional curvatures $R(x_0,y_0,y_0,x_0)$, where
$\{x_0,y_0\}$ is any orthonormal pair in $D$;

{\it mixed} sectional curvatures $R(x_0,e,e,x_0)$, where $x_0$ and
$e$ are unit vectors in $D$ and $D^{\perp}$, respectively;

{\it vertical} sectional curvature $R(e,Je,Je,e)$, where $e$ is a
unit vector in $D^{\perp}$.

Let $\{e_i\}, i=1,...,2(n-1)$ and $\{e,Je\}$ be orthonormal bases
at a point $p \in M$ of $D$ and $D^{\perp}$, respectively. Then
there arise three types of scalar curvatures associated with the
curvature tensor $R$:

{\it horizontal} scalar curvature
$$\sum_{i,j=1}^{2(n-1)} R(e_i,e_j,e_j,e_i) = \tau - 2\sigma
- 2(\sigma - \varkappa); \leqno{(2.7)}$$

{\it mixed} scalar curvature
$$\sum_{i=1}^{2(n-1)} R(e_i,e,e,e_i) = \sigma - \varkappa; \leqno{(2.8)}$$

{\it vertical} scalar curvature
$$R(e,Je,Je,e) = \varkappa. \leqno{(2.9)}$$

The invariant K\"ahler tensors
$$\pi -2\Phi +\Psi, \quad \Phi - \Psi, \quad \Psi \leqno{(2.10)}$$
are closely related to the above geometric types of sectional
curvatures because of the following their properties:

all geometric sectional curvatures of $\pi -2\Phi +\Psi$ are zero
except its horizontal sectional curvatures (more precisely this tensor
is of constant horizontal holomorphic sectional curvatures);

all geometric sectional curvatures of $\Phi - \Psi$ are zero
except its mixed sectional curvatures (more precisely this tensor is of
constant mixed sectional curvatures);

all geometric sectional curvatures of $\Psi$ are zero except its
vertical sectional curvature.

In this section we deal with holomorphic sectional curvatures
taking into account the structure of the tangent space at any point of $M$.

Any vector field $X \in {\X}M$ is decomposable in a unique way as follows:
$$X = x_0 + \tilde \eta(X)J\xi + \eta(X)\xi,$$
where $x_0$ is the projection of $X$ into ${\X}D$.

Let $span \{X, JX\}$ be an arbitrary holomorphic 2-plane (section)
in the tangent space $T_pM, \, p \in M$ generated by the unit
vector $X$. For any unit vector $Y \in span\{X,JX\}$ we have
$$\eta^2(Y) + \tilde\eta^2(Y) = \eta^2(X) + \tilde\eta^2(X).$$
The holomorphic section $span \{X, JX\}$ and the structural
section $span\{\xi,J\xi\} = D^{\perp}(p)$ form an angle
$$\varphi = \angle (span \{X, JX\}, span \{\xi, J\xi\}), \quad
\varphi \in [0, \tfrac{\pi}{2}]$$ which is uniquely determined by
the equality
$$\cos^2 \varphi = \eta^2(X) + \tilde \eta^2(X).$$
This geometric angle measures the deviation of any holomorphic
section from the structural tangent plane $D^{\perp}$.

Now we can give the basic definition in our considerations.

\begin{defn}
Let $(M,g,J,D)$ be a K\"ahler manifold with $\dim \, M = 2n \geq 4$
and $J$-invariant distribution $D$ of codimension 2. The manifold
is said to be of {\it quasi-constant holomorphic sectional
curvatures (a K\"ahler $QCH$-manifold)} if for any holomorphic section
$span \{X, JX\}$ generated by the unit tangent vector
$X \in T_pM, \, p \in M$ with
$\varphi = \angle (span \{X, JX\}, span \{\xi, J\xi \})$ the
Riemannian sectional curvature $R(X,JX,JX,X)$ may only depend on
the point $p \in M$ and the angle $\varphi$, i. e.
$$R(X,JX,JX,X) = f(p, \varphi), \quad p \in M, \,
\varphi \in [0, \tfrac{\pi}{2}].$$
\end{defn}

This notion corresponds to the notion of a Riemannian manifold of
quasi-constant sectional curvature (cf \cite{BP, GM}).

The first essential step in the study of K\"ahler $QCH$-manifolds
is to find a curvature identity characterizing these manifolds.

\begin{lem}\label{L:2.1}
Let $(M,g,J,D) \, (\dim \, M = 2n \geq 4)$ be a K\"ahler manifold of
quasi-constant holomorphic sectional curvatures. Then the
curvature tensor $R$ satisfies the following equalities
$$\begin{array}{l}
i) \, R(x_0,Jx_0,Jx_0,\xi) = R(\xi, J\xi, J\xi,x_0) = 0;\\
[2mm]
ii) \, R(x_0,Jx_0,x_0,\xi) = R(\xi, J\xi,x_0,\xi) = 0;\\[2mm]
iii) \, \displaystyle{R(x_0,J\xi,J\xi,x_0) = R(x_0,\xi,\xi,x_0) = \frac{1}{2}
R(x_0,Jx_0,J\xi,\xi)}
\end{array}$$
for all unit vectors $x_0 \in D$.
\end{lem}
{\it Proof.} Let $x_0$ be an arbitrary unit vector in $D$.

i) For any $t > 0$ we set $\displaystyle {X' = \frac{x_0 +
t\xi}{\sqrt{1 + t^2}}, \, X'' = \frac{-x_0 + t\xi}{\sqrt{1 +
t^2}}}$. The unit vectors $X'$ and $X''$ give rise to the sections
$span\{X',JX'\}$ and $span\{X'',JX''\}$ which form one and the
same angle $\varphi$ with $span\{\xi, J\xi\}$ given by the
equality
$$\cos \varphi = \frac{t}{\sqrt{1 + t^2}}.\leqno{(2.11)}$$
Hence $R(X',JX',JX',X') = R(X'',JX'',JX'',X'')$. Replacing $X'$
and $X''$ in the last equality we obtain
$$A' + A'' = A' - A'',\leqno{(2.12)}$$
where
$$A' = R(x_0,Jx_0,Jx_0,x_0) + 2t^2\{R(x_0,\xi, \xi, x_0) +
3R(x_0,J\xi,J\xi,x_0)\} + t^4R(\xi,J\xi,J\xi,\xi),$$
$$A'' = 4t\{R(x_0,Jx_0,Jx_0,\xi) + t^2R(\xi,J\xi,J\xi,x_0)\}.$$
The last equality and (2.12) imply i).

ii) In this case we consider the unit vectors
$Y' = \displaystyle {\frac{x_0 + tJ\xi}{\sqrt{1 + t^2}}}, \,
Y'' = \displaystyle{\frac{-x_0 + tJ\xi}{\sqrt{1 + t^2}}}$, $t>0.$
The holomorphic sections $span\{Y',JY'\}$ and $span\{Y'',JY''\}$ form
the same angle $\varphi$ given by (2.11) with the section
$span\{\xi,J\xi\}$. Replacing $Y'$ and $Y''$ in the equality
$R(Y',JY',JY',Y') = R(Y'',JY'',JY'',Y'')$ we find
$$B' - B'' = B' + B'', \leqno{(2.13)}$$
where
$$B' = R(x_0,Jx_0,Jx_0,x_0) + 2t^2\{R(x_0,J\xi,J\xi,x_0) +
3R(x_0,\xi,\xi,x_0)\} + t^4R(\xi,J\xi,J\xi,\xi),$$
$$B'' = 4t\{R(x_0,Jx_0,x_0,\xi) + t^2R(\xi,J\xi,x_0,\xi)\}.$$
The last equality and (2.13) imply ii).

iii) As the holomorphic sections $span\{X',JX'\}$ and
$span\{Y',JY'\}$ form the same angle $\varphi$ with the section
$span\{\xi,J\xi\}$, their sectional curvatures are equal. This
condition implies $A' = B'$ and $R(x_0,J\xi,J\xi,x_0) =
R(x_0,\xi,\xi,x_0).$ Further, the last equality in iii) follows
from the identity
$$R(x_0,Jx_0,J\xi,\xi) = R(x_0,\xi,\xi,x_0) + R(x_0, J\xi, J\xi, x_0).$$
\hfill{\bf QED}
\vskip 1mm
The curvature identity characterizing K\"ahler manifolds of
quasi-constant holomorphic sectional curvatures is given by the following

\begin{prop}\label{P:2.2}
Let $(M,g,J,D) \, (\dim \, M = 2n \geq 4)$
be a K\"ahler manifold with $J$-invariant $2(n-1)$-dimensional distribution
$D$. Then $(M,g,J,D)$ is of quasi-constant holomorphic sectional curvatures
if and only if
$$R = a\pi + b\Phi + c\Psi,$$
where $a, b$ and $c$ are functions on $M$ and the tensors $\pi,
\Phi$ and $\Psi$ are given by $(2.4), (2.5)$ and $(2.6)$,
respectively.
\end{prop}

{\it Proof.} Let $(M,g,J,D)$ be a K\"ahler manifold of
quasi-constant holomorphic sectional curvatures $f(p,\varphi)$ and
$T_pM, \, p \in M$ be the tangent space to $M$ at a fixed point $p
\in M$. Consider a unit vector $X \in T_pM$ which is neither in
$D$ nor in $D^{\perp}$, i.e. the angle $\varphi$ between $span
\{X,JX\}$ and $span\{\xi,J\xi\}$ is in the interval
$(0,\frac{\pi}{2})$. It is an easy verification that there exists
a unique unit vector $x$ in $span \{X,JX\}$ so that $x \perp \xi$
and $\tilde\eta(x) = \cos \varphi$. Then we have
$$x = \sin \varphi \, x_0 + \cos \varphi J\xi, \quad x_0 \in D, \quad
\Vert x_0 \Vert = 1 \leqno{(2.14)}$$ and
$$R(X,JX,JX,X) = R(x,Jx,Jx,x).$$
Now we replace $x$ from (2.14) into the last equality, linearize
and taking into account Lemma \ref{L:2.1} we obtain
$$\begin{array}{l}
R(X,JX,JX,X) = \\[2mm]
\sin^4 \varphi R(x_0,Jx_0,Jx_0,x_0) + 8 \sin^2 \varphi \cos^2
\varphi R(x_0,\xi,\xi,x_0) + \cos^4 \varphi R(\xi, J\xi,J\xi,\xi).
\end{array}$$
This equality can be written as
$$\begin{array}{l}
R(X,JX,JX,X) = \\[2mm]
\displaystyle{\sin^4 \varphi f\left(p,\frac{\pi}{2}\right) +
\cos^4 \varphi f(p, 0) + 8 \sin^2 \varphi \cos^2 \varphi
R(x_0,\xi,\xi,x_0)}.
\end{array}\leqno{(2.15)}$$

Further we shall express the sectional curvature
$R(x_0,\xi,\xi,x_0)$ by the functions $f(p,0), f(p,
\frac{\pi}{2})$ and $f(p,\frac{\pi}{4})$.

Let us consider the unit vector $\displaystyle{X' = \frac{x_0 +
\xi}{\sqrt{2}}}.$ We note that the angle between $span\{X',JX'\}$
and $span\{\xi,J\xi\}$ is $\displaystyle{\frac{\pi}{4}}$.
Replacing $X'$ into the equality
$$ f\left(p,\frac{\pi}{4}\right) = R(X',JX',JX',X')$$
we obtain
$$8R(x_0,\xi,\xi,x_0) = 4f\left(p,\frac{\pi}{4}\right) - f\left(p,\frac{\pi}{2}\right)
- f(p,0).$$

Now the last equality and (2.15) imply
$$\begin{array}{lc}
R(X,JX,JX,X) = & \displaystyle{f\left(p,\frac{\pi}{2}\right) +
\{4f\left(p,\frac{\pi}{4}\right) - 3f\left(p,\frac{\pi}{2}\right)
-f(p,0)\}\cos^2\varphi}\\[2mm]
& \displaystyle {+ \{2f\left(p,\frac{\pi}{2}\right) + 2f(p,0)
-4f\left(p,\frac{\pi}{4}\right)\} \cos^4 \varphi}.\end{array}$$
Putting
$$a(p) = f\left(p,\frac{\pi}{2}\right), \, b(p) = 4f\left(p,\frac{\pi}{4}\right)
- 3f\left(p,\frac{\pi}{2}\right) - f(p,0),$$
$$ c(p) = 2f\left(p,\frac{\pi}{2}\right) + 2f(p,0)-4f\left(p,\frac{\pi}{4}\right)$$
we obtain
$$R(X,JX,JX,X) = a + b\cos^2 \varphi +c\cos^4 \varphi.\leqno{(2.16)}$$
We note that the last equality is also true for the boundary
angles $\displaystyle{\varphi = 0, \frac{\pi}{2}}.$

On the other hand the holomorphic sectional curvatures of the tensor
$a\pi +b \Phi + c\Psi$ are the same as those given in (2.16). The standard
lemma for K\"ahler tensors states that if two K\"ahler tensors have the
same holomorphic sectional curvatures, they coincide (e.g. \cite {KN}).

Thus we obtained that
$$R = a \pi + b \Phi + c \Psi$$
and
$$f(p,\varphi) = a + b \cos^2\varphi + c \cos^4\varphi.$$

The inverse is an immediate verification. \hfill{\bf QED}
\vskip 2mm
Proposition \ref{P:2.2} implies in a straightforward way

\begin{cor}\label{C:2.3}
If $(M,g,J,D) \, (\dim \, M = 2n \geq 4)$ is a K\"ahler manifold of
quasi-constant holomorphic sectional curvatures, then its Ricci
tensor $\rho$ satisfies the following identity
$$\rho = \frac{\tau - 2\sigma}{2(n-1)}\,g +
\frac{2n\sigma - \tau}{2(n-1)}(\eta \otimes \eta + \tilde \eta
\otimes \tilde\eta).$$
\end{cor}

Taking into account Proposition \ref{P:2.2} and Corollary \ref{C:2.3} we
obtain expressions for the functions $a, b$ and $c$ by the geometric
functions $\varkappa, \sigma$ and $\tau$:
$$ \begin{array}{l}
\displaystyle{ a = \frac{\tau - 4 \sigma + 2\varkappa}{n(n-1)},
\quad b = \frac{4(n+2)\sigma - 2\tau - 4(n+1)\varkappa}{n(n-1)}},\\[4mm]
\displaystyle{c = \frac{\tau - 4(n+1)\sigma + (n+1)(n+2)\varkappa}{n(n-1)}}.
\end{array}\leqno{(2.17)}$$

It follows from Proposition \ref{P:2.2} that the curvature tensor $R$ of any
K\"ahler $QCH$-manifold can be represented by the tensors (2.10) in the form
$$R=a(\pi-2\Phi +\Psi)+(2a+b)(\Phi-\Psi)+(a+b+c)\Psi.\leqno{(2.18)}$$
The last formula shows that the curvature tensor $R$ of any K\"ahler
$QCH$-manifold has:

pointwise constant horizontal holomorphic sectional curvatures
$$R(x_0,Jx_0,Jx_0,x_0)=a=\frac{\tau - 2\sigma}{n(n-1)}-2\,\frac{\sigma - \varkappa}
{n(n-1)}, \quad x_0 \in D, \, \Vert x_0 \Vert = 1;\leqno{(2.19)}$$

pointwise constant mixed sectional curvatures
$$R(x_0,e,e,x_0)=\frac{2a+b}{8}=\frac{\sigma - \varkappa}{2(n-1)}, \; x_0 \in D,
\; e \in D^{\perp}, \; \Vert x_0 \Vert = \Vert e \Vert = 1;\leqno{(2.20)}$$

vertical sectional curvature
$$R(e,Je,Je,e)= a+b+c=\varkappa, \quad e\in D^{\perp},\; \Vert e \Vert = 1.\leqno{(2.21)}$$

\begin{rem}
Let $(M,g,J,D)\;(\dim M = 2n \geq 4)$be a K\"ahler manifold with $J$-invariant
distributions $D \; (\dim D = 2(n-1))$ and $D^{\perp}$. Then the structures $(g,J,D)$
generate the functions $a, b$ and $c$ by the equalities (2.17) and these functions
are related to the scalar curvatures (2.7), (2.8) and (2.9) as follows:
$$\tau - 2\sigma - 2(\sigma - \varkappa) = n(n-1)a, \leqno{(2.22)}$$
$$\sigma - \varkappa = \frac {n-1}{4}\,(2a+b), \leqno{(2.23)}$$
$$\varkappa = a + b + c. \leqno{(2.24)}$$
\end{rem}

\vskip 2mm
\section{Integrability conditions for K\"ahler manifolds of quasi-constant
         holomorphic sectional curvatures}
\vskip 2mm

The next essential step in the study of the K\"ahler $QCH$-manifolds is to
investigate the integrability conditions following from the identity
$$R = a\pi + b\Phi + c\Psi.\leqno{(3.1)}$$

First we shall make some general notes about the
structural distributions $D$ and $D^{\perp}$.

If $D^{\perp} = span \{\xi,J\xi\}$, we introduce the relative
divergences $div_0\xi$ and $div_0J\xi$ (relative codifferentials
$\delta_0\eta$ and $\delta_0\tilde\eta$) of the vector fields
$\xi$ and $J\xi$ (1-forms $\eta$ and $\tilde\eta$) with respect to
the distribution $D$:
$$div_0\xi = -\delta_0\eta = \sum_{i = 1}^{2(n-1)}(\nabla_{e_i}\eta)e_i, \quad
div_0J\xi = - \delta_0\tilde\eta = \sum_{i =
1}^{2(n-1)}(\nabla_{e_i}\tilde\eta)e_i,$$ where
$\{e_1,...,e_{2(n-1)}\}$ is an orthonormal basis of $D(p), \, p \in M$.

\begin{lem}\label{L:3.1}
Let $(M,g,J,D)$ be a K\"ahler manifold with $J$-invariant
distribution $D$ of codimension two and $D^{\perp} = span
\{\xi,J\xi\}$. Then

i) the function $(\delta_0\eta)^2 +
(\delta_0\tilde\eta)^2$ does not depend on the frame field
$\{\xi,J\xi\}$;

ii) at the points $p \in M$, where $(\delta_0\eta)^2 +
(\delta_0\tilde\eta)^2 > 0$, there exists a
geometrically determined frame field
$$\left\{\frac{\delta_0\eta \, \xi
+ \delta_0\tilde\eta \, J\xi}{\sqrt{(\delta_0\eta)^2 +
(\delta_0\tilde\eta)^2}}, \quad \frac{-\delta_0\tilde\eta \, \xi +
\delta_0\eta \, J\xi}{\sqrt{(\delta_0\eta)^2 +
(\delta_0\tilde\eta)^2}}\right\}$$ of $D^{\perp}$.
\end{lem}
{\it Proof.} Let $\xi' = \cos \varphi \, \xi + \sin \varphi \, J\xi, \quad
J\xi' = -\sin \varphi \, \xi + \cos \varphi \, J\xi$ be an
arbitrary frame field of $D^{\perp}$ and $\eta', \, \tilde\eta'$
be the 1-forms corresponding to $\xi', \, J\xi'$, respectively.
Then we find
$$\begin{array}{l}
(\nabla_{x_0}\eta')y_0 = \cos \varphi \, (\nabla_{x_0}\eta)y_0 +
\sin \varphi \,
(\nabla_{x_0}\tilde\eta)y_0,\\[2mm]
(\nabla_{x_0}\tilde\eta')y_0 = -\sin\varphi \,
(\nabla_{x_0}\eta)y_0 + \cos \varphi \,
(\nabla_{x_{0}}\tilde\eta)y_0
\end{array}$$
for all $x_0, y_0 \in D$.

The above equations immediately imply the statement of the lemma.
\hfill {\bf QED}
\vskip 1mm
We call the geometric frame field in the case ii) of the above
lemma the {\it principal} frame field for $D^{\perp}$.

In the next calculations we shall use the complexifications
$T_p^CM$ and $D^C(p), \, p\in M$ and their standard splitting
$$T_p^CM = T_p^{1,0}M \oplus T_p^{0,1}M, \quad
D^C(p) = D^{1,0}(p) \oplus D^{0,1}(p).$$ Any complex basis of $T_p^{1,0}M$
will be denoted by $\{Z_{\alpha}\}, \, \alpha = 1,...,n$ and the
conjugate basis $\{Z_{\bar\alpha} = \overline{Z_{\alpha}}\}, \,
\bar\alpha = \bar 1,...,\bar n$ will span $T_p^{0,1}M$. Counting
the special structure of the tangent spaces, we shall also
consider {\it special complex bases} of the type $\{Z_0,Z_{\lambda}\}$,
where $\displaystyle{Z_0 = \frac{\xi - iJ\xi}{2}}$ and
$\{Z_{\lambda}\}, \, \lambda = 1,...,n-1$ form a basis for
$D^{1,0}(p)$. Then $\{Z_{\bar 0}, Z_{\bar\lambda}\}, \,
\bar\lambda = \bar 1,...,\overline{n-1}$ is a special complex basis for
$T_p^{0,1}M$.

Unless otherwise stated, the Greek indices $\alpha, \beta, \gamma,
\delta, \varepsilon$ will run through $1,...,n,$ while the Greek
indices $\lambda, \mu, \nu, \varkappa, \sigma$ will run through
$1,...,n-1.$

First we give the essential components (which may not be zero) of
the fundamental tensors $g, \Omega, \eta, \tilde\eta$:
$$g_{\alpha \bar\beta} = g(Z_{\alpha},Z_{\bar\beta}), \quad
\Omega_{\alpha \bar\beta} = \Omega(Z_{\alpha},Z_{\bar\beta}) =
ig_{\alpha \bar\beta}; \leqno{(3.2)}$$
$$\eta_{\alpha} = g_{\alpha \bar 0}, \quad \eta_{\bar\alpha} =
g_{\bar\alpha 0}, \quad \tilde\eta_{\alpha} = -i\eta_{\alpha}, \quad
\tilde\eta_{\bar\alpha} = i\eta_{\bar\alpha}; \leqno{(3.3)}$$
$$\eta_{\lambda} = \eta_{\bar\lambda} = 0, \quad \eta_0 =
\eta_{\bar 0} = g_{0 \bar 0} = \frac{1}{2}. \leqno{(3.4)}$$

We introduce the following functions and 1-forms associated with
the vector fields $\nabla_{\xi}\xi$ and $\nabla_{J\xi}J\xi$:
$$p = g(\nabla_{\xi}\xi,J\xi), \quad p^* = g(\nabla_{J\xi}J\xi,\xi).
\leqno{(3.5)}$$
$$\theta(X) = g(\nabla_\xi\xi,X) - p\tilde\eta(X), \quad
\theta^*(X) = g(\nabla_{J\xi}J\xi,X) - p^*\eta(X), \quad X \in
T_pM. \leqno{(3.6)}$$ It is clear that $\theta(X) = \theta(x_0),
\, \theta^*(X) = \theta^*(x_0)$, where $x_0 = X -
\tilde\eta(X)J\xi - \eta(X)\xi.$

Taking into account that $Z_0 = \displaystyle{\frac{\xi
-iJ\xi}{2}, \, Z_{\bar 0} = \frac{\xi + iJ\xi}{2}}$ we find
$$\nabla_0\eta_{\lambda} = \frac{1}{2}(\theta_{\lambda} + \theta^*_{\lambda}),
\quad \nabla_{\bar 0}\eta_{\lambda} = \frac{1}{2}(\theta_{\lambda}
- \theta^*_{\lambda}); \leqno{(3.7)}$$
$$\nabla_0\eta_0 = \frac{p^* - ip}{4}, \quad \nabla_0\eta_{\bar 0} =
\frac{-p^* + ip}{4}.\leqno{(3.8)}$$

The next two natural statements characterize involutive
distributions $D$ and $D^{\perp}$.

\begin{lem}\label{L:3.2}
The distribution $D$ is involutive if and only if
$$\nabla_{\lambda}\eta_{\mu} - \nabla_{\mu}\eta_{\lambda} = 0, \quad
\nabla_{\lambda}\eta_{\bar\mu} = 0.$$
\end{lem}
{\it Proof}. Since the distribution $D$ is determined by the 1-forms
$\eta$ and $\tilde\eta$, it is involutive if and only if
$$d\eta(x_0,y_0) = 0, \quad d\tilde\eta(x_0,y_0) = 0, \quad x_0,y_0 \in D.$$
On the other hand, we have
$$\nabla_{\lambda}\tilde\eta_{\mu} = - i \nabla_{\lambda}\eta_{\mu}, \quad
\nabla_{\lambda}\tilde\eta_{\bar\mu} = i
\nabla_{\lambda}\eta_{\bar\mu};$$
$$d\tilde\eta_{\lambda \mu} = - i d\eta_{\lambda \mu}, \quad
d\tilde\eta_{\lambda \bar\mu} = i (\nabla_{\lambda}\eta_{\bar\mu}
+ \nabla_{\bar\mu}\eta_{\lambda}).$$ which imply the lemma.
\hfill{\bf QED}

\begin{lem}\label{L:3.3}
The distribution $D^{\perp}$ is involutive if and only if
$$\theta_{\lambda} + \theta^*_{\lambda} = 0.$$
\end{lem}
{\it Proof}. The distribution $D^{\perp}$ is involutive if and
only if the component of the Lie bracket $[\xi,J\xi]$ in $D$ is zero, i.e.
$$g(\nabla_{\xi}J\xi - \nabla_{J\xi}\xi, x_0) = 0, \quad x_0 \in D.$$
The last equality is equivalent to the condition $\theta(x_0) +
\theta^*(x_0) = 0, \, x_0 \in D.$ \hfill{\bf QED}

\begin{rem}
It is easy to check that the conditions in Lemma \ref{L:3.2} and Lemma
\ref{L:3.3} do not depend on the frame field $\{\xi,J\xi\}$.
\end{rem}
\vskip 1mm
Further we give the essential components of the tensors
$\pi, \Phi$ and $\Psi$ with respect to a complex basis
$\{Z_{\alpha}\}, \, \alpha = 1,...,n:$
$$\begin{array}{l}
\displaystyle{\pi_{\alpha\bar\beta \gamma \bar\delta} =
\frac{1}{2}(g_{\alpha\bar\beta}
g_{\gamma\bar\delta} + g_{\gamma\bar\beta}g_{\alpha\bar\delta})},\\
[2mm]\displaystyle{\Phi_{\alpha \bar\beta \gamma \bar\delta} =
\frac{1}{2}(g_{\alpha \bar\beta} \eta_{\gamma}\eta_{\bar\delta} +
g_{\gamma\bar\beta}\eta_{\alpha}\eta_{\bar\delta} + g_{\gamma
\bar\delta}\eta_{\alpha}\eta_{\bar\beta} +
g_{\alpha \bar\delta}\eta_{\gamma}\eta_{\bar\beta})},\\
[3mm]\Psi_{\alpha \bar\beta \gamma \bar\delta} = 4
\eta_{\alpha}\eta_{\gamma} \eta_{\bar\beta}\eta_{\bar\delta}.
\end{array}\leqno{(3.9)}$$

From now on in this section we assume that the manifold $(M,g,J,D)$ is
of quasi-constant holomorphic sectional curvatures, i.e. its curvature
tensor $R$ satisfies the identity (3.1). Then the essential components
of $\nabla R$ with respect to a complex basis
$\{Z_{\alpha}\},\, \alpha = 1,...,n$ are the following:
$$\begin{array}{l}
\nabla_{\alpha}R_{\beta \bar\gamma \delta \bar\varepsilon} = \\
[2mm] \displaystyle{\frac{a_{\alpha}}{2}(g_{\beta
\bar\gamma}g_{\delta \bar\varepsilon}
+ g_{\delta \bar\gamma}g_{\beta \bar\varepsilon})}\\
[2mm]
+\displaystyle{\frac{b_{\alpha}}{2}(g_{\beta\bar\gamma}\eta_{\delta}\eta_{
\bar\varepsilon} + g_{\delta \bar\gamma}\eta_{\beta}
\eta_{\bar\varepsilon} + g_{\delta
\bar\varepsilon}\eta_{\beta}\eta_{\bar\gamma}
+ g_{\beta \bar\varepsilon}\eta_{\delta}\eta_{\bar\gamma})}\\
[2mm] + \displaystyle{ \frac{b}{2}\{g_{\beta
\bar\gamma}(\nabla_{\alpha}\eta_{\delta}. \eta_{\bar\varepsilon} +
\eta_{\delta}\nabla_{\alpha}\eta_{\bar\varepsilon}) + g_{\delta
\bar\gamma}(\nabla_{\alpha}\eta_{\beta}.\eta_{\bar\varepsilon}
+ \eta_{\beta}\nabla_{\alpha}\eta_{\bar\varepsilon})}\\
[2mm] + g_{\delta
\bar\varepsilon}(\nabla_{\alpha}\eta_{\beta}.\eta_{\bar\gamma} +
\eta_{\beta}\nabla_{\alpha}\eta_{\bar\gamma}) + g_{\beta
\bar\varepsilon}(\nabla_{\alpha}\eta_{\delta}.\eta_{\bar\gamma}
+ \eta_{\delta}\nabla_{\alpha}\eta_{\bar\gamma})\}\\
[2mm]+4c_{\alpha}\eta_{\beta}\eta_{\delta}\eta_{\bar\gamma}\eta_{\bar\varepsilon} \\
[2mm] +
4c(\nabla_{\alpha}\eta_{\beta}.\eta_{\delta}\eta_{\bar\gamma}\eta_{\bar\varepsilon}
+
\eta_{\beta}\nabla_{\alpha}\eta_{\delta}.\eta_{\bar\gamma}\eta_{\bar\varepsilon}
+
\eta_{\beta}\eta_{\delta}\nabla_{\alpha}\eta_{\bar\gamma}.\eta_{\bar\varepsilon}
+
\eta_{\beta}\eta_{\delta}\eta_{\bar\gamma}\nabla_{\alpha}\eta_{\bar\varepsilon}).
\end{array} \leqno{(3.10)}$$

Now we can find the integrability conditions following from the
special form (3.1) of the curvature tensor $R$.

\begin{thm}\label{T:3.4}
Let $(M,g,J,D)$ be a K\"ahler manifold with $\dim \, M = 2n \geq 6$
and curvature tensor $R$ satisfying the identity
$$R = a\pi + b\Phi + c\Psi.$$
Then the second Bianchi identity for the tensor $R$ is equivalent
to the following conditions: \vskip 2mm i) \, $\displaystyle{da =
\frac{b \, div_0\xi}{2(n-1)} \,\eta + \frac{b \,
div_0J\xi}{2(n-1)} \, \tilde\eta;}$ \vskip 2mm ii) \,
$\displaystyle{db = \frac{(b+4c)div_0\xi}{n-1} \, \eta +
\frac{(b+4c)div_0J\xi}{n-1} \, \tilde\eta;}$ \vskip 2mm iii) \,
$b\nabla_{\lambda}\eta_{\mu} = 0, \qquad
c\nabla_{\lambda}\eta_{\mu} = 0;$ \vskip 2mm iv) \,
$\displaystyle{ b\{\nabla_{\lambda}\eta_{\bar\mu} -
\frac{div_0\xi}{2(n-1)}g_{\lambda \bar\mu} +
\frac{div_0J\xi}{2(n-1)}\Omega_{\lambda \bar\mu}\} = 0;}$ \vskip
2mm \hskip 9mm $\displaystyle{ c\{\nabla_{\lambda}\eta_{\bar\mu} -
\frac{div_0\xi}{2(n-1)}g_{\lambda \bar\mu} +
\frac{div_0J\xi}{2(n-1)}\Omega_{\lambda \bar\mu}\} = 0;}$ \vskip
2mm
 v) \, $b\theta_{\lambda} = 0, \qquad b\theta_{\lambda}^* = 0;$
\vskip 2mm vi) \, $c(\theta_{\lambda} + \theta^*_{\lambda}) =
c_{\lambda}$ \vskip 2mm \noindent with respect to a special
complex basis.
\end{thm}

{\it Proof.} First we show that the second Bianchi identity for
the curvature tensor $R = a\pi +b\Phi +c\Psi$ implies the
conditions i) - vi). The scheme of the proof is the following.
The second Bianchi identity for the curvature
tensor $R$ reduces to the equality
$$\nabla_{\alpha}R_{\beta \bar\gamma \delta \bar\varepsilon} =
  \nabla_{\beta}R_{\alpha \bar\gamma \delta \bar\varepsilon}\leqno{(3.11)}
$$
with respect to a complex basis $\{Z_{\alpha}\}, \, \alpha = 1,...,n$.

We replace $\nabla_{\alpha}R_{\beta \bar\gamma \delta
\bar\varepsilon}$ and $\nabla_{\beta}R_{\alpha \bar\gamma \delta
\bar\varepsilon}$ from (3.10) into (3.11). Further we use a special
complex basis $\{Z_0,Z_{\lambda}\}, \, \lambda = 1,...n-1,$
and substitute any of the vectors $Z_{\alpha},
\alpha = 1,...,n $ with $Z_{\lambda}, \lambda = 1,...,n-1 $ or
$Z_0$. Thus we obtain $2^5$ equalities which imply the conditions
i) - vi) of the theorem. Here we give the main points of the
proof:
\vskip 2mm
1) $\nabla_{\lambda}R_{\mu \bar\nu \varkappa \bar\sigma} =
  \nabla_{\mu}R_{\lambda \bar\nu \varkappa \bar\sigma} \quad {\To}$
$$(a_{\lambda}g_{\mu \bar\nu} - a_{\mu}g_{\lambda \bar\nu})
g_{\varkappa \bar\sigma} + (a_{\lambda}g_{\mu \bar\sigma} - a_{\mu}
g_{\lambda \bar\sigma})g_{\varkappa \bar\nu} = 0 \quad {\To}$$
$$a_{\lambda}g_{\mu \bar\nu} - a_{\mu}g_{\lambda \bar\nu} = 0
\quad {\To} \quad (n-2)a_{\lambda} = 0 \quad {\To}
\quad a_{\lambda} = 0.$$

2) $\nabla_{\lambda}R_{\mu \bar0 0 \bar\sigma} =
\nabla_{\mu}R_{\lambda \bar0 0 \bar\sigma} \quad {\To} \quad
(2a + b)_{\lambda}g_{\mu \bar\sigma} - (2a +b)_{\mu}g_{\lambda
\bar\sigma} = 0 \quad {\To}$
$$(n-2)(2a + b)_{\lambda} = 0 \quad {\To} \quad b_{\lambda}
= 0.$$

3) $\nabla_{\lambda}R_{\mu \bar\nu \varkappa \bar0} =
\nabla_{\mu}R_{\lambda \bar\nu \varkappa \bar0} \quad {\To}$
$$b(\nabla_{\lambda}\eta_{\mu} g_{\varkappa \bar\nu}
- \nabla_{\mu}\eta_{\lambda}g_{\varkappa \bar\nu} +
\nabla_{\lambda}\eta_{\varkappa}g_{\mu \bar\nu} -
\nabla_{\mu}\eta_{\varkappa}g_{\lambda \bar\nu}) = 0 \quad
{\To}$$
$$(n-2)b\nabla_{\mu}\eta_{\varkappa} = 0 \quad {\To} \quad
b\nabla_{\mu}\eta_{\varkappa} = 0.$$

4) $\nabla_{0}R_{\mu \bar0 \varkappa \bar0} = \nabla_{\mu}R_{0 \bar0
\varkappa \bar0} \quad {\To} \quad (b + 2c)\nabla_{\mu}
\eta_{\varkappa} = 0 \quad {\To} \quad
c\nabla_{\mu}\eta_{\varkappa} = 0.$
\vskip 2mm
5)$\nabla_{0}R_{\mu \bar\nu \varkappa \bar\sigma} =
\nabla_{\mu}R_{0 \bar\nu \varkappa
\bar\sigma} \quad {\To}$
$$2a_{0}(g_{\mu \bar\nu}g_{\varkappa \bar\sigma} +
g_{\varkappa \bar\nu}g_{\mu \bar\sigma}) =
b(\nabla_{\mu}\eta_{\bar\sigma}g_{\varkappa \bar\nu} +
\nabla_{\mu}\eta_{\bar\nu}g_{\varkappa \bar\sigma}) \quad
{\To}$$
$$b\nabla_{\mu}\eta_{\bar\nu}
= 2a_{0}g_{\mu \bar\nu} = \xi(a)g_{\mu \bar\nu} -
J\xi(a)\Omega_{\mu \bar\nu} \quad {\To}$$
$$\xi(a) = \frac{b\,div_0\xi}{2(n-1)}, \qquad
J\xi(a) = \frac{b\,div_0J\xi}{2(n-1)};\leqno{(3.12)}$$
$$b\{\nabla_{\mu}\eta_{\bar\nu}
- \frac{div_0\xi}{2(n-1)}g_{\mu \bar\nu} +
\frac{div_0J\xi}{2(n-1)}\Omega_{\mu \bar\nu}\} = 0.\leqno{(3.13)}$$

6) $ \nabla_{0}R_{\mu \bar0 0 \bar\sigma} = \nabla_{\mu} R_{0
\bar0 0 \bar\sigma} \quad {\To}$
$$(b + 2c)\nabla_{\mu}\eta_{\bar\sigma}
=\frac{\xi(2a + b)}{4}g_{\mu \bar\sigma} - \frac{J\xi(2a +
b)}{4}\Omega_{\mu \bar\sigma} \quad {\To}$$
$$\xi(b) = (b + 4c)\frac{div_0\xi}{n-1}, \quad
J\xi(b) = (b + 4c)\frac{div_0J\xi}{n-1}; \leqno{(3.14)}$$
$$c\{\nabla_{\mu}\eta_{\bar\sigma} - \frac{div_0\xi}{2(n-1)}g_{\mu \bar\sigma} +
\frac{div_0J\xi}{2(n-1)}\Omega_{\mu \bar\sigma}\} =
0.\leqno{(3.15)}$$

In the next three points we use the identity $(\nabla_X\eta)\xi =
0, \quad X \in {\X}M$ which implies
$$\nabla_{\alpha}\eta_0 + \nabla_{\alpha}\eta_{\bar0} = 0 \leqno{(3.16)}$$
for any $Z_{\alpha}, \alpha = 1,...,n.$ \vskip 2mm

7) $ \nabla_{0}R_{\mu \bar\nu 0 \bar\sigma} = \nabla_{\mu}R_{0
\bar\nu 0 \bar\sigma} \quad {\To}\quad b\{g_{\mu
\bar\nu}\nabla_0\eta_{\bar\sigma} + g_{\mu
\bar\sigma}\nabla_0\eta_{\bar\nu}\} = 0 \quad {\To}\quad
b\nabla_0\eta_{\bar\nu} = 0.$ \vskip 2mm \noindent Applying (3.7)
we find
$$b(\theta_{\bar\nu} - \theta^*_{\bar\nu}) = 0.$$
\vskip 2mm 8) $\nabla_{0}R_{\mu \bar0 \varkappa \bar\sigma} =
\nabla_{\mu}R_{0 \bar0 \varkappa \bar\sigma} \quad {\To} \quad
b\{g_{\varkappa \bar\sigma}\nabla_0\eta_{\mu} + g_{\mu
\bar\sigma}\nabla_0\eta_{\varkappa}\} = 0 \quad {\To} \quad
b\nabla_0\eta_{\mu} = 0.$ \vskip 2mm \noindent Applying (3.7) we
obtain
$$b(\theta_{\mu} + \theta^*_{\mu}) = 0.$$
\vskip 2mm
9) $\nabla_{0}R_{\mu \bar0 0 \bar0} = \nabla_{\mu} R_{0 \bar0 0
\bar0} \quad {\To} \quad (b + 2c)(\theta_{\mu} +
\theta^*_{\mu}) = 2 c_{\mu} \quad {\To}$
$$ c(\theta_{\mu} + \theta^*_{\mu}) = c_{\mu}.$$
Summarizing the results of 1) - 9) we obtain the equalities i) -
vi) of the theorem.

Conversely, it is easy to check that the equalities i) - vi) imply
the second Bianchi identity for $R$ (without assuming beforehand
that it is true). Thus we obtained a complete system of
integrability conditions (equivalent to the second Bianchi
identity) for the curvature tensor
$R = a\pi + b\Phi + c\Psi$.
\hfill{\bf QED}
\begin{cor}\label{C:3.5}
Let $(M,g,J,D)$ be a K\"ahler manifolds satisfying the conditions in
Theorem $\ref{T:3.4}$. If $(div_0\,\xi(p),div_0\,J\xi(p))\neq (0,0)$ at
every point $p\in M$, then the condition $b(p) \equiv 0$ implies
$c(p) \equiv 0$.
\end{cor}

Here we study K\"ahler manifolds of quasi-constant holomorphic sectional
curvatures, which have no points of constant holomorphic sectional curvatures,
i.e.
$$(b(p),c(p))\neq(0,0), \quad p \in M.\leqno{(3.17)}$$

Under this condition we have as a consequence of Theorem \ref{T:3.4}:
$$\begin{array}{l}
\displaystyle{da=\frac{b \, div_0\,\xi}{2(n-1)}\eta
+\frac{b \, div_0\,J\xi}{2(n-1)}\tilde\eta,}\\
[3mm] \displaystyle{db=\frac{(b+4c)\,div_0\,
\xi}{n-1}\,\eta+\frac{(b+4c)\,div_0\,J\xi}{n-1}\,\tilde\eta,}\\
[2mm]\nabla_{\lambda}\eta_{\mu}=0, \quad
\displaystyle{\nabla_{\lambda}\eta_{\bar\mu}=\frac{div_0\,\xi}{2(n-1)}\,
g_{\lambda\bar\mu}-\frac{div_0\,J\xi}{2(n-1)}\,\Omega_{\lambda\bar\mu},}\\
[3mm]b\theta_{\lambda} = 0, \quad b\theta^*_{\lambda}=0, \quad
c(\theta_{\lambda}+\theta^*_{\lambda})=c_{\lambda}.
\end{array}\leqno{(3.18)}$$

There arise two geometric classes of K\"ahler manifolds of quasi-constant
holomorphic sectional curvatures with respect to the distribution $D$.
\vskip 2mm
\noindent
{\bf The case:} {\it $div_0\,\xi (p)= div_0\,J\xi (p)=0$ at every point $p \in M$.}
\vskip 2mm
The equalities (3.18) reduce to the following
$$\begin{array}{l}
da = 0, \quad db = 0, \\
[2mm] \nabla_{\lambda}\eta_{\mu} = 0, \quad
\nabla_{\lambda}\eta_{\bar\mu} = 0,\\
[2mm]b\theta_{\lambda} = 0, \quad b\theta^*_{\lambda} = 0, \quad
c(\theta_{\lambda} +
\theta^*_{\lambda})=c_{\lambda}.\end{array}\leqno{(3.19)}$$

In this case the distribution $D$ is involutive according to Lemma \ref{L:3.2}
and the manifold $M$ is foliated by the integral submanifolds
of $D$.
\vskip 2mm
\noindent
{\bf The case:} {\it $(div_0 \, \xi (p), div_0 \, J\xi (p))\neq(0,0)$ at every point $p \in M.$}
\vskip 2mm
According to Corollary \ref{C:3.5}, if $b(p)\equiv 0, \, p \in M$, then $c(p) \equiv 0$,
which contradicts to (3.17). Hence $b(p)\not \equiv 0.$
\vskip 1mm
{\it The aim of the present paper is to study the class of K\"ahler
manifolds of quasi-constant holomorphic sectional curvatures
satisfying the conditions}
$$ (div_0 \,\xi (p), div_0 \, J\xi (p)) \neq (0,0), \quad b(p) \neq 0 \leqno{(3.20)} $$
at every $p \in M$.
\vskip 1mm
From (3.20) and (3.18) it follows that
$da(p) \neq 0, \, p\in M.$ Then we can choose a new
frame field $\{\xi',J\xi'\}$ for $D^{\perp}$ such that
$\displaystyle{\eta'=\frac{da}{\Vert da \Vert }}$. This frame
field satisfies the conditions $div_0 \, \xi'(p)\neq 0,div_0 \, J\xi'(p)=0,$
at any point $p \in M$.

Now the conditions (3.18) reduce to the following:
$$\begin{array}{l}
\displaystyle{da = \frac{b \, div_0 \, \xi}{2(n-1)}
\, \eta, \quad db = \frac{(b+4c)\, div_0 \, \xi}{n-1} \, \eta,}\\
[2mm] \nabla_{\lambda}\eta_{\mu}=0, \quad
\displaystyle{\nabla_{\lambda}\eta_{\bar\mu}=\frac{div_0 \,
\xi}{2(n-1)}\,g_{\lambda\bar\mu}},\\
[3mm]\theta_{\lambda} =0, \quad \theta^*_{\lambda} = 0, \quad
c_{\lambda} = 0.\end{array}\leqno{(3.21)}$$

\vskip 2mm
\section{Biconformal Transformations of K\"ahler metrics}
\vskip 2mm

It is well known that any conformal transformation
$$g'= e^{2u}g, \quad du \neq 0$$
of the metric $g$ in a K\"ahler manifold $(M,g,J)$ gives rise to a
Hermitian manifold $(M,g',J)$ which is no more K\"ahlerian .

The aim of our considerations in this section is to find the class of
$J$-invariant distributions which admit biconformal changes of the
given K\"ahler metric so that the new
metrics continue being K\"ahlerian. We obtain the
group of the biconformal transformations and
further apply it to the class of K\"ahler manifolds of
quasi-constant holomorphic sectional curvatures.

First we deal with special distributions generated by real
functions on the manifold.

Let $u$ be a real function of class $C^{\infty}$ on a K\"ahler
manifold $(M,g,J)$ and $du \neq 0$. Putting
$$\xi = \frac{grad \, u}{\Vert du \Vert}, \quad J\xi = \frac{Jgrad \, u}
{\Vert du \Vert}$$
we obtain the $J$-invariant distribution $D$, where $D^{\perp}= span \{\xi,
J\xi\}$.

With respect to a special complex basis $\{Z_0, Z_{\lambda}\},
\lambda = 1,...,n-1$ we have
$\nabla_{\lambda}u_{\bar\mu} = \Vert du \Vert
\nabla_{\lambda}\eta_{\bar\mu}$. Since
$\nabla_{\lambda}u_{\bar\mu} = \nabla_{\bar\mu}u_{\lambda}$, then
$\nabla_{\lambda}\eta_{\bar\mu}=\nabla_{\bar\mu}\eta_{\lambda}$.
The last equality implies that $\nabla_{\lambda}\eta^{\lambda} =
\nabla_{\bar\lambda}\eta^{\bar\lambda}$ which is equivalent to
$div_0J\xi = 0$.

We denote the relative Laplacian of $u$ with respect to $D$ by
$\Delta_0u$:
$$\Delta_0u = 2g^{\lambda\bar\mu}\nabla_{\lambda}u_{\bar\mu} =
\Vert du \Vert \, div_0\xi.$$

Thus any real function $u$ with $du \neq 0$ and
$\Delta_0u \neq 0$ $(div_0\xi \neq 0)$ generates a $J$-invariant
distribution $D$, where $D^{\perp} = span \{\xi, J\xi\}$, so that the
frame field $\{\xi = \frac{grad \, u}{\Vert u \Vert}, \, J\xi =
\frac{Jgrad \, u}{\Vert du \Vert}\}$ is principal. What is more,
the distribution
$$\Delta (p) = \{X \in T_pM \, \vert \, \eta (X) = 0\}, \quad p \in M$$
associated with $\eta$ (perpendicular to $\xi$) is involutive.

The above considerations give motivation to study the geometric
class of $J$-invariant distributions $D,D^{\perp} = span \{\xi,J\xi\}$
satisfying the following conditions:
$$\begin{array}{l}
i) \quad $the distribution$ \; \Delta, \; $perpendicular to$ \;
\xi, \;
$is involutive$;\\[2mm]
ii) \quad div_0 \xi \neq 0.
\end{array}\leqno{(4.1)}$$
These conditions imply that $div_0 J\xi = 0$ and the frame field $\{\xi, J\xi\}$
is principal.

Since $\Delta$ is involutive, we can consider the class of {\it
proper functions} of $\Delta$, i.e. all $C^{\infty}$-functions $q$ on $M$
satisfying the conditions (cf \cite{GM}):
$$dq \neq 0, \quad dq = \Vert dq \Vert \eta = \xi(q)\eta.\leqno{(4.2)}$$

For any proper function $q > 0$ of $\Delta$ we determine the
following transformation of the structure $(g,\eta)$:
$$\begin{array}{l}
g^* = g + (q-1)(\eta\otimes\eta + \tilde\eta\otimes\tilde\eta), \quad q > 0;\\
[2mm] \displaystyle{\eta^* = \sqrt{q} \, \eta, \quad \xi^* =
\frac{1}{\sqrt{q}} \, \xi.}
\end{array}\leqno{(4.3)}$$

The condition $q > 0$ implies $g^*$ is a Hermitian metric. We
shall call the transformations (4.3) {\it complex dilatational
transformations of the structure $(g,\eta)$}
(cf {\it dilatational transformations} in a
Riemannian manifold $(M,g,\xi)$ \cite {GM}).
We recall that the frame field $\{\xi,J\xi\}$ spanning $D^{\perp}$ is
principal and then the transformation (4.3) is an object in the
geometry of $U(n-1)\times U(1)$.

In general, the manifold $(M,g^*,J,D)$ is a Hermitian one with
metric $g^*$ and K\"ahler form $\Omega^*(X,Y) = g^*(JX,Y), \quad
X,Y \in {\X}M.$

Any complex dilatational transformation (4.3) is given by the following equalities:
$$g^*_{\beta\bar\gamma} = g_{\beta\bar\gamma} +
2(q-1)\eta_{\beta}\eta_{\bar\gamma}; \quad \eta^*_{\alpha} =
\sqrt{q} \, \eta_{\alpha}\leqno{(4.4)}$$
with respect to a complex basis $\{Z_{\alpha}\}, \alpha =1,...,n.$

Then the K\"ahler forms $\Omega$ and $\Omega^*$ of both structures
$(g,J)$ and $(g^*,J)$ are related as follows
$$\Omega^*_{\beta\bar\gamma} = \Omega_{\beta\bar\gamma} +
2i(q-1)\eta_{\beta}\eta_{\bar\gamma}.\leqno{(4.5)}$$

First we shall recall some notions and facts about the Hermitian
structures $(g^*,J)$.

The Lee form $\omega^*$ of the structure $(g^*,J)$ is determined
by
$$ \omega^*_{\alpha} =
-\frac{i}{n-1}d\Omega^*_{\alpha\beta\bar\gamma}g^{*\beta\bar\gamma}.
\leqno{(4.6)}$$

A Hermitian manifold $(M,g^*,J)$ with $\dim \, M = 2n \geq 6$ is {\it
locally conformal K\"ahler} (a {\it $W_4$-manifold} in the
classification scheme of almost Hermitian manifolds \cite {GH}) if
and only if
$$d\Omega^* = \omega^*\wedge \Omega^*.\leqno{(4.7)}$$
It is well known that (4.7) implies $\omega^*$ is closed in $\dim
M \geq 6$.

If $u$ is a local solution of the equation $2du = -\omega^*$, then
the metric $e^{2u}g^*$ is K\"ahlerian.

Now we can find the conditions for a $J$-invariant distribution $D$
of type (4.1), under which all metrics $g^*$ given by (4.3) are locally
conformal K\"ahlerian.

\begin{lem}\label{L:4.1}
Let $(M,g,J,D)$ \, $(\dim \, M = 2n \geq 6)$ be a K\"ahler
manifold with $J$-invariant distributions $D, D^{\perp} =
span\{\xi, J\xi\}$ of type $(4.1)$. Then every structure
$(g^*,J)$, where $g^*$ is a metric obtained from $g$ by a complex
dilatational transformation $(4.3)$, is locally conformal
K\"ahlerian if and only if the distributions $D, D^{\perp}$
satisfy the following conditions: \vskip 1mm i) $\theta_{\mu} +
\theta^*_{\mu} = 0;$ \vskip 1mm ii)
$\displaystyle{\nabla_{\lambda}\eta_{\bar\mu} =
\frac{div_0\xi}{2(n-1)} \, g_{\lambda\bar\mu}},$ \vskip 1mm
\noindent with respect to a special complex basis $\{Z_0,
Z_{\lambda}\}, \lambda = 1,...,n-1.$
\end{lem}

{\it Proof.} Let $g^*$ be any metric given by (4.3).

We shall use complex bases and special complex bases taking into account
the convention about the indices
$\alpha, \beta, \gamma$ and $\lambda, \mu, \nu$, respectively.

Differentiating (4.5) with respect to the Levi-Civita connection $\nabla$
of the metric $g$ and taking into account (4.2) we find
$$\begin{array}{ll}
d\Omega^*_{\alpha\beta\bar\gamma} & =
\nabla_{\alpha}\Omega^*_{\beta\bar\gamma} -
\nabla_{\beta}\Omega^*_{\alpha\bar\gamma}\\
[2mm]&= 2i(q-1)(d\eta_{\alpha\beta}\eta_{\bar\gamma} +
\eta_{\beta}\nabla_{\alpha}\eta_{\bar\gamma} -
\eta_{\alpha}\nabla_{\beta}\eta_{\bar\gamma}).
\end{array}\leqno{(4.8)}$$

The components $g^{*\beta\bar\gamma}$ of the inverse matrix of
$(g^*_{\beta\bar\gamma})$ can be found from (4.4). The
corresponding calculations give
$$g^{*\beta\bar\gamma} = g^{\beta\bar\gamma} -
\frac{2(q-1)}{q}\eta^{\beta}\eta^{\bar\gamma}.\leqno{(4.9)}$$
Replacing $d\Omega^*_{\alpha\beta\bar\gamma}$ from (4.8) and
$g^{*\beta\bar\gamma}$ from (4.9) into (4.6) we get
$$\omega^*_{\alpha} = -\frac{(q-1)div_0\xi}{n-1}\eta_{\alpha} -
\frac{q-1}{(n-1)q}(\theta_{\alpha} + \theta^*_{\alpha}).
\leqno{(4.10)}$$

Consequently
$$\omega^*_{\lambda} = - \frac{q-1}{(n-1)q}(\theta_{\lambda} +
\theta^*_{\lambda}).\leqno{(4.11)}$$

If the structure $(g^*,J)$ is locally conformal K\"ahlerian, then
$$d\Omega^*_{\alpha\beta\bar\gamma} =
\omega^*_{\alpha}\Omega^*_{\beta\bar\gamma} -
\omega^*_{\beta}\Omega^*_{\alpha\bar\gamma} =
i(\omega^*_{\alpha}g_{\beta\bar\gamma} -
\omega^*_{\beta}g_{\alpha\bar\gamma}).\leqno{(4.12)}$$
Because of (4.8) the equality (4.12) takes the form
$$2(q-1)(d\eta_{\alpha\beta}\eta_{\bar\gamma} +
\eta_{\beta}\nabla_{\alpha}\eta_{\bar\gamma} -
\eta_{\alpha}\nabla_{\beta}\eta_{\bar\gamma}) =
\omega^*_{\alpha}g^*_{\beta\bar\gamma} -
\omega^*_{\beta}g^*_{\alpha\bar\gamma}.\leqno{(4.13)}$$
Replacing
$\alpha, \beta, \gamma$ into (4.13) with $\lambda, \mu, \nu$,
respectively, we obtain
$$\omega^*_{\lambda}g^*_{\mu\bar\nu} -
\omega^*_{\mu}g^*_{\lambda\bar\nu} = 0
\quad {\To} \quad (n-2)\omega^*_{\lambda} = 0 \quad {\To} \quad
\omega^*_{\lambda} = 0.$$

Then it follows from (4.11) that
$$(q-1)(\theta_{\lambda}+\theta^*_{\lambda})=0.$$

Since the last equality is fulfilled for every proper function $q$ from (4.3),
then $\theta_{\lambda}+\theta^*_{\lambda} = 0$, which is i).

Further we replace $\alpha, \beta, \gamma$ into (4.13) with
$0, \mu, \nu$, respectively. Using (3.4) and (4.10) we get
$$-(q-1)\nabla_{\mu}\eta_{\bar\nu} = \omega^*_0g^*_{\mu\bar\nu} =
\omega^*_0g_{\mu\bar\nu} = -\frac{(q-1) \, div_0 \, \xi}{2(n-1)}\,g_{\mu\bar\nu},$$
which implies the condition ii).

Conversely, it is an immediate verification that the conditions i) and ii)
in view of (3.4) and (3.7) make (4.13) an identity.

Hence, every structure $(g^*,J)$ is locally conformal K\"ahlerian.
\hfill {\bf QED}

\begin{rem} The formula (4.10) and Lemma 3.3 imply that the following
conditions are equivalent:

i) \; $\theta_{\mu} + \theta^*_{\mu} = 0;$

ii) the distribution $D^{\perp}$ is involutive;

iii) every Lee form $\omega^*$ is collinear with the structural 1-form $\eta$.
\end{rem}

\begin{rem}
The condition ii) in the above lemma in view of the properties of
the $J$-invariant distributions $D,D^{\perp}$ is equivalent to
$$d\tilde\eta(x_0,y_0) = \frac{div_0\xi}{n-1}\Omega(x_0,y_0),
\quad x_0,y_0 \in D.\leqno{(4.14)}$$
\end{rem}

Since the function $\displaystyle{\frac{div_0\xi}{n-1} = -
\frac{\delta_0\eta}{n-1}}$ occurs quite frequently in what
follows, for the sake of brevity we shall use the denotation
$$\displaystyle{ k = \frac{div_0\xi}{n-1}}.$$

In view of Lemma \ref{L:4.1} we introduce the following notion

\begin{defn}
Let $(M,g,J,D) \, (\dim \, M = 2n \geq 6)$ be a K\"ahler manifold with
$J$-invariant distributions $D, D^{\perp} = span\{\xi,J\xi\}$.
The distribution $D$ is said to be a {\it $B$-distribution} if
the following conditions hold good:
$$\begin{array}{l}
i) \quad $the distribution$ \,\, \Delta \,\, $orthogonal to$
\,\, \xi \,\, $is involutive;$\\
[2mm]
ii) \quad $the distribution$ \,\, D^{\perp} \, $is involutive;$\\
[2mm]
iii) \quad d\tilde\eta_{\vert D} = k \,\Omega_{\vert D},
\quad k \neq 0.
\end{array}\leqno{(4.15)}$$
\end{defn}

It follows from the above definition that
$$\begin{array}{l}
\displaystyle{\nabla_{\lambda}\eta_{\mu} =
\nabla_{\mu}\eta_{\lambda}, \quad
\nabla_{\lambda}\eta_{\bar\mu} = \frac{k}{2}\, g_{\lambda\bar\mu}},\\
[2mm] \displaystyle{\frac{div_0\xi}{n-1}=k \neq 0, \quad div_0J\xi
= 0.}
\end{array}$$

Thus we have:

\vskip 2mm {\it The geometry of a $B$-distribution in a K\"ahler
manifold is determined by: \vskip 1mm the symmetric tensor
$\nabla_{\lambda}\eta_{\mu}$; \vskip 1mm the 1-form $\theta$;
\vskip 1mm the functions $k, p$ and $p^*$.} \vskip 2mm

Now we can introduce transformations of the tensors $g, \eta$ changing
the K\"ahler structure $(g,J)$ into a K\"ahler structure $(g',J)$.

Let $(M,g,J,D)\,(\dim \, M = 2n \geq 6)$ be a K\"ahler manifold with
$B$-distribution $D \; (D^{\perp} = span \{\xi,J\xi\})$. If
$$g^* = g + (q-1)(\eta \otimes \eta + \tilde\eta \otimes
\tilde\eta), \quad q > 0, \quad dq = \xi(q)\eta$$
is a complex dilatational change of the metric $g$, then $(g^*,J)$ is a locally
conformal K\"ahler structure whose Lee form $\omega^*$ in view of (4.10) is given by
$$\omega^* = -(q-1)k\eta.$$
If $u$ is a local solution of the equation
$$2du = - \omega^* = k(q-1)\eta, $$
then $u$ is a proper function of the distribution $\Delta$ and
$2\xi(u) = k(q-1)$. Putting
$$g' = e^{2u}g^* = e^{2u}\{g + (q-1)(\eta\otimes\eta +
\tilde\eta\otimes\tilde\eta)\}$$
we obtain $g'$ is a K\"ahler metric.

It is clear that the metric $g'$ is determined by the proper function $q$
up to a constant factor.

Since $q > 0$ we set $q = e^{2v}$ and obtain another form of the metric $g'$:
$$g' = e^{2u}\{g + (e^{2v}-1)(\eta \otimes \eta + \tilde \eta \otimes
\tilde \eta)\}.$$

\begin{defn}
Let $(M,g,J,D) \, (\dim \, M = 2n \geq 6)$ be a K\"ahler manifold
with $B$-distribution $D$. The transformation
$$\begin{array}{l}
g' = e^{2u}\{g - (\eta\otimes\eta + \tilde\eta\otimes\tilde\eta)\}
+ e^{2(u + v)}(\eta\otimes\eta +
\tilde\eta\otimes\tilde\eta);\\
[2mm]
\eta' = e^{u+v}\eta, \quad \xi' = e^{-(u+v)}\xi; \\
[2mm]dv = \xi(v)\eta \neq 0, \quad 2du = k (e^{2v} - 1)\eta.
\end{array}\leqno{(4.16)}$$
is said to be a {\it biconformal transformation} of the structure $(g,\eta)$.
\end{defn}
From the above definition it follows that $(M,g',J,D)$ is a K\"ahler manifold
with the same distributions $D, D^{\perp}$ and $\Delta$. By a given proper
function $v$ the metric $g'$ is determined up to a constant factor (homothety).

\begin{prop}\label{P:4.2}
Let $(M,g,J,D) \, (\dim \, M = 2n \geq 6)$ be a K\"ahler
manifold with $B$-distribution $D$ and $(g', \eta')$ be any structure
obtained from the structure $(g, \eta)$ by a biconformal transformation.
Then $D$ continues being a $B$-distribution with
respect to the new structure $(g',\eta')$.
\end{prop}

{\it Proof.} Let the structure $(g', \eta')$ be given by (4.16). Since
the distributions $\Delta$ and $D^{\perp}$ do not change, then
the first two conditions of Definition 4.4 are fulfilled.

To prove the third condition of Definition 4.4 we differentiate the equality
$\tilde\eta' = e^{u+v}\tilde\eta$ and taking into account (4.16) we get
$$d\tilde\eta'_{\vert D}= e^{u+v}d\tilde\eta_{\vert D} = e^{u+v}k\Omega_{\vert D}
= e^{v-u}k \Omega'_{\vert D}.$$
Hence
$$d\tilde\eta'_{\vert D} = k'\Omega'_{\vert D}, $$
where
$$k'=e^{v-u}k. \leqno{(4.17)}$$
\hfill{\bf QED}

\begin{prop}\label{P:4.3}
Let $(M,g,J,D) \, (\dim \, M = 2n \geq 6)$ be a K\"ahler
manifold with $B$-distribution $D$. Then the biconformal
transformations of the structure $(g,\eta)$ form a group.
\end{prop}
{\it Proof.} I. Let $v$ be a proper function of the distribution
$\Delta$, i.e. $dv = \xi (v)\eta \neq 0$ and
$$g' = e^{2u}\{g - (\eta\otimes\eta +
\tilde\eta\otimes\tilde\eta)\} + e^{2(u+v)}(\eta\otimes\eta +
\tilde\eta\otimes\tilde\eta),$$
$$ \eta' = e^{u+v}\eta, \quad
du = \frac{k(e^{2v}-1)}{2} \, \eta.$$
By direct computations taking into account (4.17) we find the inverse
transformation
$$g = e^{-2u}\{g' - (\eta'\otimes\eta' +
\tilde\eta'\otimes\tilde\eta')\} + e^{-2(u+v)}(\eta'\otimes\eta' +
\tilde\eta'\otimes\tilde\eta'),$$
$$\eta = e^{-(u+v)}\eta', \quad d(-u) = \frac{k'(e^{-2v}-1)}{2} \, \eta'.$$

II. Let $v'$ and $v''$ be proper functions of the distribution $\Delta$
generating the following biconformal transformations:
$$g' = e^{2u'}\{g - (\eta\otimes\eta +
\tilde\eta\otimes\tilde\eta)\} + e^{2(u'+v')}(\eta\otimes\eta +
\tilde\eta\otimes\tilde\eta),$$
$$\eta' = e^{u'+v'}\eta, \quad du'= \frac{k(e^{2v'}-1)}{2} \, \eta;$$
$$g'' =e^{2u''}\{g' - (\eta'\otimes\eta' +
\tilde\eta'\otimes\tilde\eta')\} +
e^{2(u''+v'')}(\eta'\otimes\eta' +
\tilde\eta'\otimes\tilde\eta'),$$
$$\eta''= e^{u''+v''}\eta', \quad du'' = \frac{k'(e^{2v''}-1)}{2} \, \eta'.$$
It is easy to check that
$$g'' = e^{2(u'+u'')}\{g- (\eta\otimes\eta +
\tilde\eta\otimes\tilde\eta)\} +
e^{2(u'+u''+v'+v'')}(\eta\otimes\eta +
\tilde\eta\otimes\tilde\eta),$$
$$\eta''= e^{u'+u''+ v'+v''}\eta, \quad d(u'+u'') = \frac{k(e^{2(v'+v'')}-1)}{2} \, \eta.$$
Taking into account I and II we conclude that the biconformal
transformations (4.16) form a group, which is generated by the additive group
of the real $C^{\infty}$-functions on $M$. \hfill{\bf QED}

\vskip 2mm
\section{Tensor invariants of the biconformal group of transformations}
\vskip 2mm

Let $(M,g,J,D)$ be a K\"ahler manifold with $B$-distribution $D$. In order to find
tensor invariants of the group of biconformal transformations we shall further
specialize the distribution $D$.

Taking into account that the $J$-invariant distribution of any K\"ahler
$QCH$-manifold of the class under consideration is a $B$-distribution satisfying
the conditions $\nabla _{\lambda}\eta_{\mu} = 0$ and $\theta = \theta^* = 0$
we introduce the following type of $J$-invariant distributions:
\begin{defn}
Let $(M,g,J,D)$ be a K\"ahler manifold with $J$-invariant distribution
$D$ of codimension two. The distribution $D$ is said to be a {\it $B_0$-distribution}
if it satisfies the following conditions:
$$\begin{array}{l}
1) \; D \; $is a $B$-distribution;$\\
[2mm]
2) \; (\nabla_{x_0}\eta)(y_0)- (\nabla_{Jx_0}\eta)(Jy_0) = 0, \quad x_0, y_0 \in D;\\
[2mm] 3) \; \theta = 0.
\end{array}\leqno{(5.1)}$$
\end{defn}

From the above definition it follows immediately that $\theta^* = 0$.

Under the assumptions (5.1) we shall find additional conditions
for the basic functions $p$ and $p^*$, introduced by (3.5), and
$\displaystyle{k = \frac{div_0 \, \xi}{n-1}}.$

First we precise the form of the tensor $\nabla \eta$.

\begin{lem}\label{L:5.1}
Let $(M,g,J,D)$ be a K\"ahler manifold with $B_0$-distribution $D$. Then
$$\begin{array}{ll}
(\nabla_X\eta)Y =& \displaystyle{\frac{k}{2}\{g(X,Y)
- \eta(X)\eta(Y) - \tilde\eta(X)\tilde\eta(Y)\}}\\
[2mm]& + p \eta(X)\tilde\eta(Y) - p^* \tilde\eta(X)\tilde\eta(Y);
\quad X,Y \in {\X}M.\end{array}$$
\end{lem}

{\it Proof.} We consider a special complex basis
$\{Z_0, Z_{\lambda}\}, \lambda = 1,...,n-1$. Then the assertion of the lemma
is equivalent to the equalities:
$$\begin{array}{l}
\nabla_{\lambda}\eta_{\mu} = 0, \quad
\displaystyle{\nabla_{\lambda}\eta_{\bar\mu}= \frac{k}{2} \,
g_{\lambda\bar\mu}};\\
[2mm] \nabla_0\eta_{\mu}=0, \quad \nabla_{\bar 0}\eta_{\mu} = 0;
\end{array}\leqno{(5.2)}$$
$$\begin{array}{l}
\nabla_{\lambda}\eta_0 = 0, \quad \nabla_{\lambda}\eta_{\bar 0}= 0;\\
[2mm] \displaystyle{\nabla_0\eta_0 = -\nabla_0\eta_{\bar 0} =
\frac{1}{4}(p^*-ip)}.
\end{array}\leqno{(5.3)}$$

The first two equalities of (5.2) are given by the conditions 1) and 2) of (5.1).
The last two equalities of (5.2) follow from the formulas (3.7) and the
condition 3).

Since the distribution $\Delta$ orthogonal to $\xi$ is involutive, then
$d\eta(x_0,J\xi) = 0$ for any $x_0 \in D$. The last equality and the condition
$\theta^* = 0$ imply $g(\nabla_{x_0}\xi, J\xi)=0, \; x_0 \in D$. From here we
obtain the first two equalities of (5.3). The remaining equalities of (5.3)
are given by the formulas (3.8).
\hfill {\bf QED}
\vskip 2mm
In the next lemma we obtain integrability conditions concerning the functions
$k, \; p$ and $p^*$ as a consequence of the conditions
$\nabla_{\lambda}\eta_{\mu} = 0$ and $\theta = \theta^* = 0$.

\begin{lem}\label{L:5.2}
Let $(M,g,J,D) \; (\dim \, M = 2n \geq 6)$ be a K\"ahler manifold with
$B_0$-distribution. Then the following equalities
$$p = 0, \quad dk = -k(k+p^*)\eta, \quad dp^* = \xi(p^*)\eta$$
hold good.
\end{lem}
{\it Proof.} Applying Lemma \ref{L:5.1} we obtain
$$d\eta = p\eta \wedge \tilde\eta;\leqno{(5.4)}$$
$$d\tilde\eta=k\Omega - (k+p^*)\eta\wedge\tilde\eta. \leqno{(5.5)}$$
The exterior derivation of (5.4) because of (5.5) gives
$$k p(\eta_{\alpha}g_{\beta\bar\gamma}-\eta_{\beta}g_{\alpha\bar\gamma})
=2(p_{\alpha}\eta_{\beta}-p_{\beta}\eta_{\alpha})\eta_{\bar\gamma}$$
with respect to a complex basis $\{Z_{\alpha}\}, \alpha = 1,...,n.$

Passing on a special complex basis $\{Z_0, Z_{\lambda}\}, \lambda = 1,...,n-1$
we substitute $\alpha = 0, \; \beta = \mu, \; \gamma = \nu$ in the above
equality and find $kp = 0$. Since $k \neq 0$, then
$$p = 0.$$
A similar operation applied to (5.5) gives as a result the equality
$$\begin{array}{l}
\{k_{\alpha}+k(k+p^*)\eta_{\alpha}\}g_{\beta\bar\gamma}-
\{k_{\beta}+k(k+p^*)\eta_{\beta}\}g_{\alpha\bar\gamma}\\
[2mm]-2\{(k+p^*)_{\alpha}\eta_{\beta}-
(k+p^*)_{\beta}\eta_{\alpha}\}\eta_{\bar\gamma} = 0.\end{array}
\leqno{(5.6)}$$

Proceeding to a special complex basis in (5.6) we substitute
$\alpha = \lambda, \; \beta = \mu, \; \gamma = \nu$ and get
$$k_{\lambda}g_{\mu\bar\nu}-k_{\mu}g_{\lambda\bar\nu} = 0 \quad
{\To} \quad (n-2)k_{\lambda} = 0 \quad {\To} \quad k_{\lambda} = 0.$$

The substitution $\alpha = 0, \; \beta = \mu, \; \gamma = \nu$
into (5.6) gives as a consequence
$$\{k_0+\frac{1}{2}k(k+p^*)\}g_{\mu\bar\nu}=0,$$
which implies $J\xi(k)=0, \quad \xi(k)=-k(k+p^*).$

Thus we obtained the equalities
$$k_{\lambda} = 0, \quad J\xi(k)=0, \quad \xi(k)=-k(k+p^*),\leqno{(5.7)}$$
which are equivalent to
$$dk = \xi(k)\eta = -k(k+p^*)\eta.$$

The exterior derivation in the last equality gives $dp^*\wedge \eta = 0$, i.e.
$$dp^*=\xi(p^*)\eta.$$
\hfill{\bf QED}

Now, taking into account Lemma \ref{L:5.2}, we have from Lemma \ref{L:5.1}
$$\nabla \eta =
\frac{k}{2}(g-\eta\otimes\eta-\tilde\eta\otimes\tilde\eta)
-p^*\tilde\eta\otimes\tilde\eta,\leqno{(5.8)}$$
where the function $p^*$ in view of (5.7) satisfies the relation
$$p^* = -\frac{\xi(k)+k^2}{k}.\leqno{(5.9)}$$

The equality (5.8) is equivalent to the equalities
$$\nabla_{\alpha}\eta_{\beta} = p^*\eta_{\alpha}\eta_{\beta},
\leqno{(5.10)}$$
$$\nabla_{\alpha}\eta_{\bar\beta}=\frac{k}{2}g_{\alpha\bar\beta}
-(k+p^*)\eta_{\alpha}\eta_{\bar\beta}\leqno{(5.11)}$$
with respect to a complex basis $\{Z_{\alpha}\}, \alpha = 1,...,n.$
\begin{lem}\label{L:5.3}
Let $(M,g,J,D) \; (\dim \, M=2n \geq 6)$ be a K\"ahler manifold with
$B_0$-distribution $D$. Then the following equalities hold good:
$$\sigma= \frac{1}{2k}\xi(k^2+2kp^*)+\frac{n+1}{2}(k^2+2kp^*);
\leqno{(5.12)}$$
$$\varkappa = \frac{1}{2k}\xi(k^2+2kp^*)+k^2+2kp^*;\leqno{(5.13)}$$
$$R(X,Y)\xi = 4\frac{\sigma-\varkappa}{n-1}\{\Phi(X,Y)\xi-
\Psi(X,Y)\xi\}+\varkappa\Psi(X,Y)\xi.\leqno{(5.14)}$$
\end{lem}

{\it Proof.} Writing the equality (5.8) in the form
$$\nabla_X\xi=\frac{k}{2}\{X-\eta(X)\xi-\tilde\eta(X)J\xi\}
-p^*\tilde\eta(X)J\xi, \quad X\in {\X}M,$$
we get in a straightforward way
$$\begin{array}{ll}
R(X,Y)\xi = & \displaystyle{-\frac{1}{4}(k^2+2kp^*)\{\eta(X)Y-\eta(Y)X}\\
[3mm]
& \displaystyle{-\tilde\eta(X)JY+\tilde\eta(Y)JX+2g(JX,Y)J\xi\}}\\
[2mm] & \displaystyle{-\frac{1}{2k}\xi(k^2+2kp^*)
\{\eta(X)\tilde\eta(Y)-\eta(Y)\tilde\eta(X)\}J\xi},
\end{array}\leqno{(5.15)}$$
$X,Y \in {\X}M$.

Taking a trace in (5.15) we find
$$\rho(X,\xi)=\{\frac{1}{2k}\xi(k^2+2kp^*)+\frac{n+1}{2}(k^2+2kp^*)\}\eta(X),
\quad X \in {\X}M,\leqno{(5.16)}$$
which implies (5.12).

Further we replace $X=J\xi, Y=\xi$ into (5.15) and find (5.13).

Taking into account (5.15), (2.5), (2.6), (5.12) and (5.13) we obtain (5.14).
\hfill {\bf QED}

\begin{rem}
Let $D$ be a $B_0$-distribution. Then the equality (5.16) means that $\xi$
is an eigen vector field for the Ricci operator $\rho$.

Further the equality (5.14) implies that the manifold under consideration
is of pointwise constant mixed sectional curvatures
$$R(x_0,\xi,\xi,x_0)= \frac{\sigma-\varkappa}{2(n-1)}=\frac{1}{4}(k^2+2kp^*),
\quad x_0 \in D,\; \Vert x_0 \Vert = 1. \leqno{(5.17)}$$

What is more, the equalities (5.17), (5.13) and (5.9) express that
the function $k$ completely determines the mixed sectional
curvatures and the vertical sectional curvature of the manifold.
\end{rem}

The next assertion justifies the introduction of the class of $B_0$-distributions.
\begin{prop}\label{P:5.4}
Let $(M,g,J,D) \; (\dim \, M = 2n \geq 6)$ be a K\"ahler manifold with
$B_0$-distribution $D$ and $D^{\perp}=span\{\xi,J\xi\}$. If
$(g',\eta')$ is any structure obtained from the given structure
$(g,\eta)$ by a biconformal transformation, then
$D$ is also a $B_0$-distribution with respect to the new structure
$(g',\eta')$.
\end{prop}
{\it Proof.} Let
$$\begin{array}{c}
g'=e^{2u}\{g +
(e^{2v}-1)\,(\eta\otimes\eta+\tilde\eta\otimes\tilde\eta)\},
\quad \eta'=e^{u+v}\eta;\\
[2mm] \displaystyle{dv=\xi(v)\,\eta, \quad
du=\frac{k(e^{2v}-1)}{2}\,\eta}
\end{array}\leqno{(5.18)}$$
be a biconformal change of the structure $(g,\eta)$ and $\nabla'$ be the
Levi-Civita connection of the metric $g'$.

According to Proposition \ref{P:4.2} the distribution $D$ is a $B$-distribution with
respect to the structure $(g',\eta')$. It remains to prove that $\theta' = 0$ and
$(\nabla'_{x_0}\eta')(y_0)-(\nabla'_{Jx_0}\eta')(Jy_0) = 0, \quad x_0, y_0 \in D.$

By direct computations we find the relation between $\nabla'$ and $\nabla$:
$$\begin{array}{ll}
\nabla'_XY=&\nabla_XY + \xi(u)\{\eta(X)Y + \eta(Y)X +
\tilde\eta(X)JY + \tilde\eta(Y)JX\}\\
[2mm]
& + \xi(v-u)\{(\eta(X)\eta(Y) - \tilde\eta(X)\tilde\eta(Y))\xi \\
[2mm]& + (\eta(X)\tilde\eta(Y)+\tilde\eta(X)\eta(Y))J\xi\}; \quad
X,Y \in {\X}M. \end{array}\leqno{(5.19)}$$

Then we find the components of $\nabla'\eta'$ with respect to a complex basis:
$$\nabla'_{\alpha}\eta'_{\beta} =
e^{-(u+v)}\{p^*-\xi(u+v)\}\,\eta'_{\alpha}\eta'_{\beta},\leqno{(5.20)}$$
$$\nabla'_{\alpha}\eta'_{\bar\beta}=e^{v-u}\frac{k}{2}
(g'_{\alpha\bar\beta}-2\eta'_{\alpha}\eta'_{\bar\beta})
-e^{-(u+v)}\{p^*-\xi(u+v)\}\,\eta'_{\alpha}\eta'_{\bar\beta}.\leqno{(5.21)}$$

Now let $\{Z'_{\lambda}, Z'_0\}, \lambda = 1,..., n-1$ be a special complex basis
with respect to the structure $(g', \eta')$. Taking into account (5.20), (5.21) and
(3.7) we get
$$\nabla'_{\lambda}\eta'_{\mu} = 0, \quad \nabla'_{\lambda}\eta'_0 =
\frac{\theta'_{\lambda}+\theta'^{*}_{\lambda}}{2}=0, \quad
\nabla'_{\lambda}\eta'_{\bar0}=\frac{\theta'_{\lambda}-\theta'^{*}_{\lambda}}{2}=0,$$
which implies the assertion. \hfill{\bf QED}
\vskip 2mm
The equality (5.21) implies the change (4.17) of the function $k$ and
the change of the function $p^*$:
$${p^*}' =e^{-(u+v)}\{p^*-\xi(u+v)\}.\leqno{(5.22)}$$

Now we can prove the main theorem in this section.
\begin{thm}\label{T:5.5}
Let $(M,g,J,D) \; (\dim \, M = 2n \geq 6)$ be a K\"ahler manifold with
$B_0$-distribution $D$ and $D^{\perp}=span\{\xi,J\xi\}$. Then the
tensor of type $(1.3)$
$$QC(R) = R-a\pi-b\Phi-c\Psi$$
is a biconformal invariant.
\end{thm}
{\it Proof.} Let $(g',\eta')$ be a structure obtained from the given structure $(g,\eta)$
by the biconformal transformation (5.18). Writing the operator $QC(R)$ of type (1,3)
in the form
$$QC(R)= R-a(\pi - 2\Phi + \Psi)-(2a+b)(\Phi - \Psi)-(a+b+c)\Psi \leqno{(5.23)}$$
we shall prove
$$QC(R')= QC(R).\leqno{(5.24)}$$

Taking into account (5.19) we find the relation between the curvature tensors
$R'$ and $R$ of type (1.3):
$$\begin{array}{ll}
R'-R=&-2k\xi(u)\pi -4k\xi(v-u)\Phi_1-4\{\xi^2(u)-(k+p^*)\xi(u)\}\Phi_2\\
[2mm]
&-\{\xi^2(v-u)-(2k+p^*)\xi(v-u)\}\Psi,\end{array}\leqno{(5.25)}$$
where
$$\begin{array}{l}
\Phi_1(X,Y)Z=\\
[2mm]
\displaystyle{\frac{1}{8}\{g(Y,Z)(\eta(X)\xi+\tilde\eta(X)J\xi)
-g(X,Z)(\eta(Y)\xi+\tilde\eta(Y)J\xi)}\\
[2mm]
+g(JY,Z)(\eta(X)J\xi-\tilde\eta(X)\xi)-g(JX,Z)(\eta(Y)J\xi-\tilde\eta(Y)\xi)\\
[2mm] -2g(JX,Y)(\eta(Z)J\xi-\tilde\eta(Z)\xi)\};\end{array}$$
$$\begin{array}{l}
\Phi_2(X,Y)Z=\\
[2mm]
\displaystyle{\frac{1}{8}\{(\eta(Y)\eta(Z)+\tilde\eta(Y)\tilde\eta(Z))X
-(\eta(X)\eta(Z)+\tilde\eta(X)\tilde\eta(Z))Y}\\
[2mm] +(\eta(Y)\tilde\eta(Z)-\tilde\eta(Y)\eta(Z))JX
-(\eta(X)\tilde\eta(Z)-\tilde\eta(X)\eta(Z))JY\\
[2mm]
-2(\eta(X)\tilde\eta(Y)-\tilde\eta(X)\eta(Y))JZ\},\end{array}$$
$X,Y,Z \in {\X}M.$

It follows that $\Phi_1+\Phi_2=\Phi$.

The equality (5.25) can be rewritten in the form
$$\begin{array}{ll}
R'-R=&-2k\xi(u)(\pi - 2\Phi +\Psi)\\
[2mm] &\displaystyle{-4k\xi(v)(\Phi_1-\frac{1}{2}\Psi)
-4(\xi^2(u)-p^*\xi(u))(\Phi_2-\frac{1}{2}\Psi)}\\
[3mm] &-(\xi^2(u+v)-p^*\xi(u+v))\Psi.\end{array}\leqno{(5.26)}$$

It is easy to check that the tensors
\;$\displaystyle{\pi - 2\Phi +\Psi, \; \Phi_1-\frac{1}{2}\Psi,\;
\Phi_2-\frac{1}{2}\Psi}$ and $\Psi$ of type (1,3) change in the following way:
$$\pi'- 2\Phi' + \Psi'=e^{2u}(\pi - 2\Phi + \Psi),$$
$$\Phi'_1-\frac{1}{2}\Psi'=e^{2u}(\Phi_1-\frac{1}{2}\Psi),$$
$$\Phi'_2-\frac{1}{2}\Psi'=e^{2(u+v)}(\Phi_2-\frac{1}{2}\Psi),$$
$$\Psi'=e^{2(u+v)}\Psi.$$

To prove the equality (5.24) because of (5.26) it is sufficient to prove
the following equalities
$$-2k\xi(u)= e^{2u}a'-a, \leqno{(5.27)}$$
$$-4k\xi(v)= e^{2u}(2a'+b')-(2a+b), \leqno{(5.28)}$$
$$-4(\xi^2(u)-p^*\xi(u))=e^{2(u+v)}(2a'+b')-(2a+b), \leqno{(5.29)}$$
$$-(\xi^2(u+v)-p^*\xi(u+v))=e^{2(u+v)}(a'+b'+c')-(a+b+c). \leqno{(5.30)}$$

Considering the equality (5.25) onto the distribution $D$ we have
$$R'(x_0,y_0)z_0=R(x_0,y_0)z_0-2k\xi(u)\pi(x_0,y_0)z_0, \quad x_0,y_0,z_0 \in D.$$
Taking traces in the last equality twice we find
$$e^{2u}\{\tau'-2\sigma'-2(\sigma'-\varkappa')\}=\tau-2\sigma-2(\sigma-\varkappa)
-2n(n-1)k\xi(u).$$
This equality gives (5.27) because of (2.22).

Further we apply the formula (5.17) and taking into account (4.17), (5.22) we find
$$e^{2u}\frac{\sigma'-\varkappa'}{n-1}=\frac{\sigma-\varkappa}{n-1}-k\xi(v).$$
The last equality becomes (5.28) because of (2.23).

From Lemma \ref{L:5.3} it follows that
$$e^{2(u+v)}\sigma'= \sigma -\{\xi^2(nu+v)-p^*\xi(nu+v)\},$$
$$e^{2(u+v)}\varkappa'= \varkappa - \{\xi^2(u+v)-p^*\xi(u+v)\}.\leqno{(5.31)}$$
Hence
$$e^{2(u+v)}\frac{\sigma'-\varkappa'}{n-1}=\frac{\sigma-\varkappa}{n-1}-
\{\xi^2(u)-p^*\xi(u)\},$$
which gives (5.29).

Finally the equality (5.31) in view of (2.24) becomes (5.30).
\hfill{\bf QED}
\vskip 2mm
Theorem \ref{T:5.5} implies in a straightforward way that the Ricci trace
$\rho (QC(R))$ is a biconformal invariant of type (0, 2).
\begin{cor}\label{C:5.6}
Let $(M,g,J,D) \; (\dim \, M = 2n \geq 6)$ be a K\"ahler manifold with
$B_0$-distribution $D$. Then the tensor of type $(0,2)$
$$\rho-\frac{\tau-2\sigma}{2(n-1)}\,g-\frac{2n\sigma-\tau}{2(n-1)}
(\eta\otimes\eta+\tilde\eta\otimes\tilde\eta)$$
is a biconformal invariant.
\end{cor}

The following relative scalar invariant is important to the next
applications of Theorem \ref{T:5.5}.

\begin{cor}\label{C:5.7}
Let $(M,g,J,D) \; (\dim \, M = 2n \geq 6)$ be a K\"ahler manifold with
$B_0$-distribution $D$ and $D^{\perp}=span\{\xi,J\xi\}$. If
$(g',\eta')$ is obtained from $(g,\eta)$ by the biconformal
transformation $(5.18)$, then
$$a' + k'^2 = e^{-2u}(a + k ^2).$$
\end{cor}

{\it Proof.} From (5.27) and (5.18) we have
$$e^{2u}a'=a-2k\xi(u)=a-k^2(e^{2v}-1).$$
On the other hand (4.17) gives
$$e^{2u}k'^2=e^{2v}k^2,$$
which implies the assertion.
\hfill {\bf QED}
\vskip 2mm
Thus there arise three classes of $B_0$-distributions, which are
invariant under the biconformal group of transformations. These
classes are determined by the conditions
$$a+k^2>0, \quad a+k^2=0, \quad  a+k^2<0,\leqno{(5.32)}$$
respectively.

The next question to clear up is to describe $B_0$-distributions in
the standard flat K\"ahler manifold ${\C}^n$.

Let $(M,g',J,D) \; (\dim \, M = 2n \geq 6)$ be a K\"ahler
manifold with flat Levi-Chivita connection $\nabla'$ of $g'$ and
$B_0$-distribution $D$ \; ($D^{\perp}=span\{\xi',J\xi'\}$).
It follows from (5.8) that
$$ d\eta'=0;\leqno{(5.33)}$$
$$\nabla'_{\xi'}\xi'=0.\leqno{(5.34)}$$
The condition $R'=0$ and the equalities (5.17), (5.9) imply
$$\xi'(k')+\frac{k'^2}{2}=0.\leqno{(5.35)}$$
Then from (5.8) and (5.35) we obtain
$$\nabla'_x\xi'=\frac{k'}{2}x, \quad x\in \Delta.\leqno{(5.36)}$$

Under the conditions (5.35) and (5.36) we can describe
$B_0$-distributions in ${\C}^n=\{Z=(z^1,...,z^n)\}$.

Let $({\C}^n,g',J,D) \; (n \geq 3)$ be the complex Euclidean space
endowed with the standard flat K\"ahler structure $(g',J)$. If $D$
is a $B_0$-distribution in ${\C}^n$ , then the condition (5.33)
allows us to put locally $\eta'=dt$. In this case (5.35) becomes
$\displaystyle{\frac{dk'}{dt}=-\frac{k'^2}{2}}$ whose general
solution is $\displaystyle{k'=\frac{2}{t+t_0}, \; t_0 = const.}$
Let us consider any integral submanifold $S^{2n-1}$ of the
integrable distribution $\Delta$. Since $\xi'$ is a unit normal
field to $S^{2n-1}$, the condition (5.36) implies $S^{2n-1}$ is
(part of) a hypersphere with radius $\displaystyle{r=\frac{2}{\vert
k' \vert}= \vert t+t_0 \vert}.$ On the other hand according to
(5.34) the integral curves of the vector field $\xi'$ are straight
lines, which pass trough the center $Z_0$ of $S^{2n-1}$. Hence
$S^{2n-1}$ are concentric hyperspheres.

Choosing $Z_0$ as the origin $O$ of ${\C}^n$ we obtain the following
\vskip 2mm
{\it Canonical example of a flat K\"ahler manifold with
$B_0$-distribution:}
\vskip 2mm
$M'= {\C}^n \setminus \{O\}; \; (g', J)$ is the standard flat
K\"ahler structure in ${\C}^n$; \; $\displaystyle{\xi' =
\frac{Z}{\Vert Z \Vert}}$, where $Z$ is the position vector of the
corresponding point in $M'$.
\vskip 2mm
It is easy to check that the distribution $D$ determined by
$D^{\perp}= span\{\xi',J\xi'\}$ is a $B_0$-distribution and
$$\eta' = dr, \quad k'=\frac{2}{r}; \quad
r^2=\sum_{\alpha =1}^nz^{\alpha}z^{\bar\alpha}, \quad r > 0.$$

\begin{defn}
Let $(M,g,J,D) \; (\dim \, M = 2n \geq 6)$ be a K\"ahler manifold with
$B_0$-distribution $D$ and $D^{\perp}=span\{\xi,J\xi\}$. The
manifold is said to be {\it biconformally flat} if there
exists a flat metric $g'$ obtained from the metric $g$ by a
biconformal transformation. The metric $g$ is said to be
{\it a biconformally flat} metric.
\end{defn}
\begin{thm}\label{T:5.8}
Let $(M,g,J,D) \; (\dim \, M = 2n \geq 6)$ be a K\"ahler manifold with
$B_0$-distribution $D$ and $D^{\perp}=span\{\xi,J\xi\}$. The
manifold is biconformally flat if and only if the
following conditions hold good
$$QC(R) = 0, \quad a+k^2>0.$$
\end{thm}

{\it Proof.} Let the structure $(g',\eta')$ be obtained from the
given structure $(g,\eta)$ by the biconformal transformation (5.18).
If $R'=0$, then $a'=b'=c'=0$ and $QC(R')=0$. According to
Theorem \ref{T:5.5} \, $QC(R)=0$.

On the other hand Corollary \ref{C:5.7} gives that
$$a + k^2 = e^{2u}k'^2>0.$$

For the inverse, let $QC(R)=0$ and $a+k^2>0$. We shall construct a
new structure $(g',\eta')$ of type (5.18), whose curvature tensor $R'=0$.

From Corollary \ref{C:5.7} and (4.17) we have
$$e^{2u}a'= a-(e^{2v}-1)k^2.$$
Then the condition $a'=0$ is equivalent to
$$2v=\ln\frac{a+k^2}{k^2}.\leqno{(5.37)}$$
Lemma \ref{L:5.2} and (3.21) imply that the function $v$ in (5.37) satisfies the condition
$dv=\xi(v)\eta$.

According to (5.8) $d\eta=0$ and we have locally $\eta = ds$
for some function $s$. Replacing $v$ from (5.37) and $\eta = ds$ into the equality
$$2du=k(e^{2v}-1)\eta$$
we find
$$2u=\int\frac{a}{k} \; ds. \leqno{(5.38)}$$

Let us consider the biconformal transformation (5.18) determined by the functions
(5.37) and (5.38). According to Theorem \ref{T:5.5} and (5.37) we have
$$R'=b'\Phi'+c'\Psi'.$$
Taking into account (2.23), (2.24) and $a'=0$ the conditions $b'=c'=0$ are
equivalent to $\sigma'=\kappa'=0$. According to (5.12) and (5.13) the last
equalities are equivalent to the equality $k'^2+2k'{p^*}'=0.$

From (5.28) because of (2.23) and (5.17) we have
$$e^{2u}(k'^2+2k'{p^*}')=k^2+2kp^*-2k\xi(v).$$
Calculating $2k\xi(v)$ from (5.37) we obtain
$$e^{2u}(k'^2+2k'{p^*}')=-\frac{k}{a+k^2}\{\xi(a)+ak-k(k^2+2kp^*)\}.$$
The right hand side of the last equality is identically zero because of
the formulas
$$\xi(a)=\frac{kb}{2}, \quad 2a+b=2(k^2+2kp^*),$$
which follow from (3.21), (2.23) and (5.17).

Thus we proved that $a'=b'=c'=0$. Hence $R'=0$.
\hfill{\bf QED}
\vskip 2mm
Theorem \ref{T:5.8} allows to obtain locally all biconformally flat K\"ahler metrics.

Let $(M'={\C}^n\setminus\{O\},g',J,D)\;(n\geq3)$ be the canonical flat K\"ahler
manifold with $B_0$-distribution $D$. Then
$$r^2=2g'_{\alpha\bar\beta}z^{\alpha}z^{\bar\beta},\quad
g'_{\alpha\bar\beta}=\left \{
\begin{array}{l}
\frac{1}{2} \quad \alpha = \beta \\[1mm]
0 \quad \alpha \neq \beta,
\end{array}\right. \quad \eta'=dr, \quad k'= \frac{2}{r}.$$

For an arbitrary function $v(r^2)\in C^{\infty}$ we construct the metric
$$\begin{array}{c}
g=e^{-2u}\,\{g'+(e^{-2v}-1)(\eta'\otimes\eta'+\tilde\eta'\otimes\tilde\eta')\};\\
[2mm]
\displaystyle{\frac{d(-u)}{dr^2}=\frac{e^{-2v}-1}{2r^2},}\end{array}\leqno{(5.39)}$$
which is determined up to a constant factor (homothety).

Theorem \ref{T:5.8} implies
\begin{cor}\label{C:5.9}
Let $(M,g,J,D) \; (\dim M = 2n\geq6)$ be a K\"ahler manifold with $B_0$-distribution $D$.
If the metric $g$ is biconformally flat, then it can be presented locally
in the form $(5.39)$ and vice versa.
\end{cor}

Finally we show that the metrics (5.39) are closely related to the K\"ahler metrics $g$
in ${\C}^n$ whose potential function is $f(r^2) \in C^{\infty}$, i.e.
$$g=\partial\bar\partial f(r^2).\leqno{(5.40)}$$
These metrics have been used as a source for K\"ahler metrics with special properties
(cf \cite{TL, B}).

\begin{thm}\label{T:5.10}
Let $(M,g,J,D) \; (\dim M = 2n\geq6)$ be a K\"ahler manifold with $B_0$-distribution $D$.
Then the metric $g$ is biconformally flat if and only if it is of type $(5.40)$.
\end{thm}

{\it Proof.} Let the metric $g$ be given by (5.39). Putting
$$f(r^2)=\frac{1}{2}\int{e^{-2u}}\,dr^2$$
we obtain $g$ is of type (5.40).

For the inverse, let the function $f(r^2)$ generate the K\"ahler metric $(5.40)$.
Putting $\displaystyle{e^{-2u}=2f'>0, \quad e^{-2v}=1+r^2\frac{f''}{f'}>0,}$
we obtain the metric $g$ is of type (5.39). An easy check shows that
$\displaystyle{\xi'(-u)=\frac {k'}{2}(e^{-2v}-1).}$ Hence the metric $g$ is
biconformally flat.
\hfill {\bf QED}

\vskip 2mm
\section{K\"ahler structures on rotational hypersurfaces}
\vskip 2mm
The aim of this section is to show that any rotational hypersurface
$(M^{2n},\bar g)$ which has no common points with the axis of revolution
carries a geometrically determined complex structure $J$ and $(M^{2n},\bar g,J)$
can be considered as a locally conformal K\"ahler manifold. Further we show that
$M^{2n}$ carries a natural K\"ahler metric $g$, so that $(M^{2n},g,J)$ is a
K\"ahler manifold of quasi-constant holomorphic sectional curvatures.

Let $Oe$ be a fixed coordinate system in ${\R}$ and $M^{2n}$ be
a rotational hypersurface in ${\R}^{2n+1}={\C}^n\times {\R}$ with axis of
revolution $l={\R}$. We consider the class of {\it rotational hypersurfaces
having no common points with the axis of revolution.} Then $M^{2n}$ is a
one-parameter family of spheres $S^{2n-1}(s),\,s\in I \subset {\R}$ considered
as hyperspheres in ${\C}^n$ with corresponding centers $q(s)e$ on $l$ and
radii $t(s)>0$. If $Z$ is the radius vector of any point $p \in M^{2n}$ with
respect to the origin $O$, then the unit outer normal $n$ of the parallel
$S^{2n-1}(s)$ at the point $p$ is
$$n=\frac{Z-q(s)e}{t(s)}.$$
Hence
$$Z=q(s)e+t(s)n.\leqno{(6.1)}$$

Further we assume that $s$ is a natural parameter for the meridian
$$\gamma(s)=q(s)e+t(s)n \leqno{(6.2)}$$
in the plane $Oen$ ($n$ - fixed), i. e. $q'^2+t'^2=1$.

Because of (6.2) and (6.1) the unit tangent vector field $\bar\xi$ to
the meridian $\gamma(s)$ is
$$\bar \xi=\frac{d\gamma}{ds}=q'e+t'n=\frac{\partial Z}{\partial s}.
\leqno{(6.3)}$$

Since the normal to $M^{2n}$ lies in the plane $Oen$, we always choose
the unit vector field $N$ normal to $M^{2n}$ by the condition that the
couples $(e,n)$ and $(\bar\xi,N)$ have the same orientation. Then
taking into account (6.3), we have
$$N=-t'e+q'n.$$

Let $\bar g$ be the standard metric on ${\R}^{2n+1}={\C}^n \times {\R}$
with flat Levi-Civita connection $\nabla'$. We denote the induced metric
on $M^{2n}$ by the same letter $\bar g$. Then $(M^{2n},\bar g,\bar\xi)$
becomes the Riemannian warped product manifold \cite {BO, GM}:
$$\gamma(s)\times_{t(s)}S_0^{2n-1},$$
where $S_0^{2n-1}$ is the unit hypersphere in ${\C}^n$, centered at the
origin $O$.

Denoting by $\bar\nabla$ the Levi-Civita connection on $(M^{2n},\bar g)$
we have \cite {GM}:
\vskip 1mm
$$\begin{array}{l}
\vspace{2mm}
\bar\nabla_x\bar\xi=\displaystyle{\frac{t'}{t}\,x}; \quad x\in {\X}M^{2n}, \,\,
\bar g(x,\bar\xi)=0,\\
\vspace{2mm}
\bar\nabla_{\bar\xi}\bar\xi=0.\end{array}
\leqno{(6.4)}$$

Let $\bar\eta$ be the 1-form corresponding to the unit vector field $\bar\xi$
with respect to the metric $\bar g$, i. e. $\bar\eta(X)=\bar g(\bar\xi,X),
\,\, X\in {\X}M^{2n}$. If $\bar \pi$ and $\bar \Phi$ are the tensors

$$\begin{array}{ll}
\vspace{2mm}
\bar\pi (X,Y)Z= & \bar g(Y,Z)X-\bar g(X,Z)Y,\\
\vspace{2mm}
\bar \Phi (X,Y)Z= & \bar g(Y,Z)\bar\eta(X)\bar\xi-\bar g(X,Z)\bar\eta(Y)\bar\xi\\
\vspace{1mm}
& +\bar\eta(Y)\bar\eta(Z)X-\bar\eta(X)\bar\eta(Z)Y, \quad X,Y,Z \in {\X}M^{2n},
\end{array}$$
then the curvature tensor $\bar R$ of the rotational hypersurface $M^{2n}$ has the
form \cite{GM}:
$$\bar R=\frac{1-t'^2}{t^2}\,\bar\pi-\frac{1-t'^2+tt''}{t^2}\,\bar\Phi.\leqno{(6.5)}$$
This equality implies that the rotational hypersurface $M^{2n}$ is
conformally flat.

We shall introduce a complex structure on any rotational hypersurface $M^{2n}$
having no common points with the axis $l$.

First we consider the almost contact Riemannian structure on the parallels of
the rotational hypersurface $M^{2n}$ induced from the corresponding ${\C}^n$.

Let $(J_0,\bar g)$ be the standard flat K\"ahler structure of any ${\C}^n$
considered as a hyperplane in ${\R}^{2n+1}$ perpendicular to the axis $l$.
Then any parallel $S^{2n-1}(s)$ ($s$ - fixed) being a hypersphere in
$({\C}^n,J_0,\bar g)$ carries a natural almost contact Riemannian structure
$(\varphi, \tilde{\bar\xi}, \tilde{\bar\eta},\bar g)$ determined as follows
\cite{T1,T2}:
$$\begin{array}{l}
\vspace{2mm}
\tilde{\bar\xi}=J_0n;\\
\vspace{2mm}
\tilde{\bar\eta}(x)=\bar g(x,\tilde{\bar\xi}),
\quad x\in{\X}S^{2n-1}(s);\\
\vspace{2mm}
\varphi(x)=J_0x+\tilde{\bar\eta}(x)n.
\end{array}\leqno{(6.6)}$$

The corresponding Weingarten and Gauss formulas of the imbedding
$S^{2n-1}(s)\subset{\C}^n$ are
$$\begin{array}{l}
\vspace{2mm}
\nabla'_xn=\displaystyle{\frac{1}{t}\,x};\\
\vspace{2mm}
\nabla'_xy=\tilde\nabla_xy-\displaystyle{\frac{1}{t}\,\bar g(x,y)n};
\quad x,y \in {\X}S^{2n-1}(s),\end{array} \leqno{(6.7)}$$
where $\tilde\nabla$ is the induced Levi-Civita connection on the sphere
$S^{2n-1}(s)$.
From (6.6) and (6.7) it follows directly that
$$\tilde\nabla_x \tilde{\bar\xi}=\frac{1}{t}\,\varphi x,\quad x \in{\X}S^{2n-1}(s);
\leqno{(6.8)}$$
$$\begin{array}{ll}
\vspace{2mm}
0=(\nabla'_xJ_0)y= & \displaystyle{(\tilde\nabla_x \varphi)y-
\frac{1}{t}\,\tilde{\bar\eta}(y)x+\frac{1}{t}\,\bar g(x,y)\tilde{\bar\xi}}\\
\vspace{2mm}
& \displaystyle{-\left(\frac{1}{t}\,\bar g(x,\varphi y)+
(\tilde\nabla_x \tilde{\bar\eta})y\right)n},
\quad x,y \in {\X}S^{2n-1}(s).\end{array}$$
Hence
$$(\tilde\nabla_x \varphi)y=\frac{1}{t}\,\left(\tilde{\bar\eta}(y)x-\bar g(x,y)
\tilde{\bar\xi}\right), \quad x,y \in {\X}S^{2n-1}(s).\leqno{(6.9)}$$

If the structure $(\varphi,\tilde{\bar\xi},\tilde{\bar\eta},\bar g)$ of an almost
contact Riemannian manifold satisfies the conditions
\vskip 1mm
$$\begin{array}{l}
\vspace{2mm}
\tilde\nabla_x\tilde{\bar\xi}=\alpha \,\varphi x, \quad \alpha=const;\\
\vspace{1mm}
(\tilde\nabla_x\varphi)y=\alpha \left(\tilde{\bar\eta}(y)x-
\bar g(x,y)\tilde{\bar\xi}\right),\end{array}$$
then the manifold is called an {\it $\alpha$-Sasakian manifold} \cite{JV}.
In the case $\alpha = 1$ these manifolds are the usual Sasakian manifolds.

Taking into account the equalities (6.8) and (6.9) we conclude that the
structure (6.6) on any hypersphere $S^{2n-1}(s)$ ($s$ - fixed) in ${\C}^n$ is
$\displaystyle{\frac{1}{t}}$-Sasakian ($\displaystyle{\frac{1}{t}=\alpha}
=const$).

Now we can introduce a complex structure $J$ on the rotational hypersurface
$M^{2n}$ subordinated to the orientation $\bar\xi$ of the meridians.

Let $T_pM^{2n}$ be the tangent space to $M^{2n}$ at any point $p$. Then the
vector fields $\bar\xi$ and $\tilde{\bar\xi}$ defined by (6.6) determine a
distribution $D$ so that $D^{\perp}=span\{\bar\xi,\tilde{\bar\xi}\}$. We define
the almost complex structure $J$ associated with $\bar\xi$ as follows:
$$J_{|D}=J_0, \quad J\bar\xi=\tilde{\bar\xi}, \quad J\tilde{\bar\xi}=-\bar\xi.
\leqno{(6.10)}$$
It is clear that $J$ is an almost complex structure and $(M^{2n},\bar g,J)$ becomes
an almost Hermitian manifold.

The almost complex structure $J$ defined by (6.10) and the structure $\varphi$
on $S^{2n-1}(s)$ given by (6.6) are related as follows:
$$\varphi(x)=Jx+\tilde{\bar\eta}(x)\bar\xi, \quad x \in {\X}S^{2n-1}(s).
\leqno{(6.11)}$$

Below we give the Weingarten and Gauss formulas for the embedding
$S^{2n-1}(s)\subset M^{2n}$ with normal vector field $\bar\xi$:
$$\begin{array}{l}
\vspace{2mm}
\displaystyle{\bar\nabla_x\bar\xi=\frac{t'}{t}\,x},\\
\vspace{2mm}
\bar\nabla_xy=\tilde\nabla_xy-\displaystyle{\frac{t'}{t}\,\bar g(x,y)\bar\xi},
\quad x,y \in {\X}S^{2n-1}(s).
\end{array}\leqno{(6.12)}$$

The Gauss formula and (6.8) imply immediately that
$$\bar\nabla_x \tilde{\bar\xi}=\frac{1}{t}\,\varphi x-\frac{t'}{t}\tilde{\bar\eta}(x)
\bar\xi, \quad x \in {\X}S^{2n-1}(s).\leqno{(6.13)}$$

Further we prove
\begin{prop}
Let $(M,\bar g)$ be a rotational hypersurface in ${\C}^n\times R$ with axis of
revolution $l={\R}$, which has no common points with $l$ and the meridians of
$M^{2n}$ are oriented with the unit vector field $\bar\xi$ . If $J$ is the
almost complex structure $(6.10)$ associated with $\bar\xi$, then the covariant
derivative of $J$ satisfies the identity
$$(\bar\nabla_XJ)Y=\frac{t'-1}{t}\left(\bar g(X,Y)\tilde{\bar\xi}-
\tilde{\bar\eta}(Y)X-\bar\eta(Y)JX+\bar g(JX,Y)\bar\xi\right) \leqno{(6.14)}$$
for all vector fields $X,Y\in{\X}M^{2n}$.
\end{prop}

{\it Proof}: To prove (6.14) it is sufficient to show that the following
equalities hold good:
$$\begin{array}{l}
\vspace{2mm}
(i) \, (\bar\nabla_{\bar\xi}J)x=0, \quad \bar g(\bar\xi,x)=0;\\
\vspace{2mm}
(ii) \, (\bar\nabla_{\bar\xi}J)\bar\xi=0;\\
\vspace{2mm}
(iii) \, (\bar\nabla_xJ)y=\displaystyle{\frac{t'-1}{t}\left(\bar g(x,y)
\tilde{\bar\xi}-\tilde{\bar\eta}(y)x+\bar g(\varphi x,y)\bar\xi\right)},
\quad \bar g(\bar\xi,x)=\bar g(\bar\xi,y)=0;\\
\vspace{1mm}
(iv) \, (\bar\nabla_xJ)\bar\xi=\displaystyle{\frac{1-t'}{t}
\left(Jx+\tilde{\bar\eta}(x)\bar\xi \right)},
\quad \bar g(\bar\xi,x)=0.
\end{array}$$

The equality (i) follows because of the fact that the hyperplanes
${\C}^n(s)$ are parallel along any meridian $\gamma(s)$.

Taking into account (i) and (6.4)it follows that (ii) holds good.

The equality (iii) follows from (6.12), (6.13), (6.11) and (6.9).

Finally (iv) follows from (6.12), (6.13) and (6.11). \hfill{\bf QED}
\vskip 2mm
The identity (6.14) implies that the almost complex structure $J$
is integrable, i. e. $(M^{2n},\bar g,J)$ is a Hermitian manifold.
Moreover, this manifold is in the class $W_4$ according to the
classification in \cite {GH}. A simple computation shows that the
Lee form of the manifold is $\displaystyle{\frac{1-t'}{t}\,\bar\eta}$.
From (6.4) it follows that the 1-form $\bar\eta$ is closed and therefore
the Lee form of $(M^{2n},\bar g,J)$ is also closed. Then the manifold
under consideration is locally conformal K\"ahler in all dimensions
$2n \geq 4$. This implies $(M^{2n},\bar g,J)$ carries a conformal
K\"ahler metric which is flat by virtue of the fact that $\bar g$ is
conformally flat.

Our aim is to define another nontrivial K\"ahler metric on
$(M^{2n},\bar g,J)$, which is naturally determined by its geometric
structures.

We constrain the class of the rotational hypersurfaces which have no
common points with the axis of revolution, assuming the inequality
$t'(s)\neq 0, \, s \in I$. This condition means that the meridian
$\gamma(s)$ of $M^{2n}$ has no points in which the tangents are
parallel to the axis $l$. Under the condition $t'\neq 0$ we can choose
in a unique way the orientation $\bar\xi$ of the meridians so that
$t'(s)>0$.

In what follows we consider the class of rotational hypersurfaces
$(M^{2n},\bar g,J)$ satisfying the following inequalities:
$$t(s)>0, \quad t'(s)>0; \quad s\in I. \leqno{(6.15)}$$

Under the conditions (6.15) we construct the metric
$$g=\bar g+(t'-1)(\bar\eta\otimes\bar\eta+\tilde{\bar\eta}
\otimes\tilde{\bar\eta}). \leqno{(6.16)}$$

Taking into account (6.8) and (6.11) it is easy to check that the
K\"ahler form of the metric $g$ is closed, i. e. $g$ is a K\"ahler
metric. We call $g$ a {\it complex dilatational K\"ahler} metric on the
rotational hypersurface $M^{2n}$. This metric is relevant to the subject
of our considerations in the present paper because of the following
\begin{thm}\label{T:6.2}
Let $(M^{2n},\bar g,J,\bar\xi)\,(2n\geq 4)$ be a rotational hypersurface
satisfying the conditions $(6.15)$. Then the complex dilatational
K\"ahler metric $g$ given by $(6.16)$ is of quasi-constant holomorphic
sectional curvatures.
\end{thm}

{\it Proof:} Let $\nabla$ be the Levi-Civita connection of the metric
(6.16). Using (6.4) and (6.13) we find in a standard way that
\vskip 1mm
$$\begin{array}{ll}
\vspace{2mm}
\nabla_XY = & \bar\nabla_XY + \displaystyle{\frac{1-t'}{t}}
\left\{\bar\eta(JX)JY+\bar\eta(JY)JX-\bar g(X,Y)\bar\xi\right\}\\
\vspace{2mm}
& \displaystyle{+\left(\frac{1-t'}{t}+\bar\xi(\ln{\sqrt {t'}})\right)}
\left\{\bar\eta(X)\bar\eta(Y)-\bar\eta(JX)\bar\eta(JY)\right\} \bar\xi\\
\vspace{2mm}
& \displaystyle{-\left(\frac{1-t'}{t}+\bar\xi(\ln{\sqrt {t'}})\right)}
\left\{\bar\eta(X)\bar\eta(JY)+\bar\eta(JX)\bar\eta(Y)\right\}J\bar\xi
\end{array}\leqno{(6.17)}$$
for all $X,Y \in {\X}M^{2n}$.

Taking into account (6.17), (6.5) and (6.16) we compute the curvature tensor
$R$ of the connection $\nabla$:
$$R=a\pi+b\Phi+c\Psi,$$
where
$$a=\frac{4(1-t')}{t^2}, \quad b=8\left(\frac{t'-1}{t^2}-\frac{t''}{2tt'}\right),
\quad c=\frac{4(1-t')}{t^2}+\frac{5t''}{2tt'}+\frac{t''^2-t't'''}{2t'^3}.
\leqno{(6.18)}$$

Applying Proposition \ref {P:2.2} we obtain the assertion.
\hfill{\bf QED}
\vskip 2mm
Since $t'^2=1-q'^2$, then $t'\in (0,1]$. Hence $a\geq 0$ in (6.18).

From the equality $q'^2=1-t'^2$ it also follows that the function $t=t(s)$
determines the rotational hypersurface $M^{2n}$ up to a translation along
the axis $l$ and a symmetry with respect to the hyperplane ${\C}^n$ through
the origin $O$.

Let $p\in M^{2n}$ and $\gamma (p), \, S^{2n-1}(p)$ be the corresponding meridian
and parallel through the point $p$. From (6.14), (6.5) and (6.18) it follows that
the following conditions are equivalent:
$$\begin{array}{l}
\vspace{2mm}
(1) \, \, $the tangent$ \, \, \bar\xi \, \, $to$ \, \, \gamma \, \, $at$ \, \,
p \in M^{2n} \, \, $is perpendicular to the axis$ \, \, l;\\
\vspace{2mm}
(2) \,\, g=\bar g \, \, $on the parallel$ \, \, S^{2n-1}(p) \, \,
$through the point$ \,\, p;\\
\vspace{2mm}
(3) \,\, \bar\nabla J=0 \, \, $on$ \, \, S^{2n-1}(p);\\
\vspace{2mm}
(4) \,\, \bar R=0 \, \, $on$ \, \, S^{2n-1}(p);\\
\vspace{2mm}
(5) \,\, R=0 \, \, $on$ \, \, S^{2n-1}(p).\end{array}$$

As a consequence of Theorem \ref{T:6.2} we can find the rotational hypersurfaces
$M^{2n}$ whose complex dilatational K\"ahler metric is of constant holomorphic
sectional curvatures.

Let $b=0$ in (6.18). Then Corollary \ref{C:3.5} implies that $c=0$ and the
metric $g$ is of constant holomorphic sectional curvature $a=const >0$.

Solving the equation
$$b=8\left(\frac{t'-1}{t^2}-\frac{t''}{2tt'}\right)=0$$
we obtain the meridian in the form $q=q(t)$.

Considering the meridian in the usual coordinate system $Oxy$ with axis of
revolution $l=Oy$ we have:
\begin{prop}
Any rotational hypersurface $M^{2n}$ which carries a complex dilatational
K\"ahler metric of constant holomorphic sectional curvature $a=const >0$
is generated by a meridian of the type
$$\gamma : y=\pm\frac{1}{\sqrt a}\left(\sqrt{8-ax^2}+\ln{\frac{\sqrt{8-ax^2}-
2}{\sqrt{8-ax^2}+2}}\right)+y_0, \quad 0<x<\frac{2}{\sqrt a}.$$
\end{prop}

\end{document}